\documentclass[11pt]{article}

\usepackage{amsmath}
\usepackage{amscd}
\usepackage{latexsym}
\usepackage[all]{xy}
\newfont{\gothique}{eufm10 scaled 1100}  %% gotic
\newcommand{\goth}[1]{\mbox{\gothique{#1}}}

\newcommand{\PP}{{\bf{P}}}
\newcommand{\CC}{\bf{C}}
\newcommand{\AL}{\alpha}
\newcommand{\ARR}{\longrightarrow}
\newcommand{\LK}{L+K_X}
\newcommand{\LKS}{\mid L+K_X \mid}
\newcommand{\ZE}{Z_e}

\newcommand{\ZD}{{\cal{Z}}}
\newcommand{\ZDO}{\stackrel{\circ}{\cal{Z}}}
\newcommand{\EE}{{\cal E}}
\newcommand{\SEE}{{\EE},[e]}
\newcommand{\ZA}{([Z],[\alpha])}
\newcommand{\JA}{{\bf{J}}( X;L,d)} 
\newcommand{\JAA}{{\bf{J}}}

\newcommand{\XD}{X^{[d]}}
\newcommand{\EXT}{Ext^{n-1}_{Z_e}(L)}
\newcommand{\JUE}{{\bf{J}}_{U_{\cal E}}}
\newcommand{\JUEO}{{\bf{J}}_{\RUE}}
\newcommand{\EXX}{Ext^{n-1}}
\newcommand{\OO}{{\cal O}}
\newcommand{\OL}{{\OO}_X (L)}
\newcommand{\UE}{U_{\cal E}}
\newcommand{\UEP}{\UE^{\prime}}
\newcommand{\RU}{\stackrel{\circ}{U}}
\newcommand{\RUE}{{\stackrel{\circ}{U}}_{\EE}}
\newcommand{\RUEP}{{\stackrel{\circ}{U}}^{\prime}_{\EE}}
\newcommand{\EXTS}{{\bf{Ext}}^{n-1}_{U_{\cal E}}}
\newcommand{\EXTSU}{{\bf{Ext}}^{n-1}_{U}}
\newcommand{\EXTSO}{{\bf{Ext}}^{n-1}_{\RU_{\cal E}}}
\newcommand{\EEF}{{\cal E}{\it{xt}}}
\newcommand{\EX}{Ext^{n-1}} 
\newcommand{\JB}{{\bf{J}}} 
\newcommand{\EA}{([e],[\alpha])}
\newcommand{\HT}{{\bf \tilde{H}}}
\newcommand{\HB}{{\bf{H}}}
\newcommand{\ID}{{\cal{I}}}
\newcommand{\RELTAN}{{\cal T}_{\JBB / \RUE}}
\newcommand{\PIPO}{\pi^{\ast}_{\RUE}\big( p_{{\RUE}{\ast}} \OO_{\ZDO}\big)}
\newcommand{\PIP}
{ \pi^{\ast}_{\UE}\left(p_{\UE \ast}\big(\OO_{\ZD}\big) \right)}
\newcommand{\PIPU}
{ \pi^{\ast}_{U}\left(p_{U \ast}\big(\OO_{\ZD_U}\big) \right)}
\newcommand{{\JBB}}{\breve{\JB}}
\newcommand{{\JBBU}}{{\breve{\JB}}_{\RUE}}

\newcommand{\QB}{{\bf q}}

\newcommand{\HOZ}{H^0(\OO_{Z_e})}
\newcommand{\HOZK}{H^0(\OO_{Z_e} (\LK))}
\newtheorem{thm}{Theorem}[section]
\newtheorem{lem}[thm]{Lemma}
\newtheorem{pro}[thm]{Proposition}
\newtheorem{cor}[thm]{Corollary}
\newtheorem{rem}[thm]{Remark}

\newtheorem{defi}[thm]{Definition}  %[section]

\newlength{\myskip}
\newenvironment{pf}{
     \addvspace{\myskip}  

     \noindent {\it Proof.$\, $}}
     {$\Box$

     \addvspace{\myskip}
     }

\makeatletter
\renewcommand{\@seccntformat}[1]{\S \/ {\csname the#1\endcsname}\hspace{0.5em}}
\makeatother

\title{CONFIGURATIONS OF POINTS AND STRINGS}
\author{Igor Reider}

\begin{document}
\bibliographystyle{amsplain}

\maketitle

\setcounter{section}{-1}
\numberwithin{equation}{section}

\begin{abstract}
Let $X$ be a smooth projective variety of dimension $n \geq2$.
It is shown that a finite configuration of points on $X$ subject to
certain geometric conditions possesses rich inner structure.
On the mathematical level this inner structure is a variation of
Hodge-like structure. 
 As a consequence one can attach to such point configurations:
\\
\indent
(i) Lie algebras and their representations
\\
\indent
(ii) Fano toric variety whose hyperplane sections are Calabi-Yau varieties.
\\
These features imply that the points cease to be 0-dimensional objects and acquire dynamics
of linear operators ``propagating" along the paths of a particular trivalent graph. Furthermore following particular linear operators along the ``shortest" paths of the graph
one creates, for every point of the configuration, a distinguished hyperplane section of the Fano variety in (ii), i.e. the points ``open up" to become Calabi-Yau varieties.
Thus one is led to a picture which is very suggestive of quantum gravity according to string theory.
\end{abstract}

\section{Introduction} 
The theory of algebraic curves or Riemann surfaces is one of the most beautiful and harmonious chapters of Algebraic geometry. One of the reasons for this is undoubtedly due to the object called the Jacobian of a curve. One knows that it is an abelian variety naturally associated with a smooth projective curve. Virtually everything one wants to know about a curve can be extracted from its Jacobian.
The idea to use the Jacobian to study correspondences of algebraic curves goes back to A.Weil, \cite{[W]}.
It is with this insight in mind that we tried to develop the theory of Jacobian for algebraic surfaces. In \cite{[R1]}, motivated by the
problem of correspondences of algebraic surfaces, we introduced an object called nonabelian Jacobian of a smooth complex projective surface.
Our approach can be viewed as a natural generalization of the classical theory. Namely, viewing the Jacobian of a curve as a parameter space
for line bundles with a fixed Chern class one is tempted to generalize this for an $n$-dimensional ($n\geq2$) smooth projective variety along
the following lines:
\\
\\
\indent
fix the Chern classes for a holomorphic vector bundle of rank $n$ over a complex projective variety $X$ of dimension $n$ and construct
a ``canonical" family of such bundles. The parameter space of such a family might be considered as a nonabelian version of Jacobian for
$n$-dimensional varieties.
\\
\\
The paper \cite{[R1]} treats the case of smooth projective surfaces, i.e. $n=2$. 
To explain the goals and contents of the present paper we begin by a brief
 summary of the results in {\it loc. cit}. 
\\
\indent 
Let  $X$ be a smooth complex projective surface. For 
a fixed divisor  $L$  (up to the rational equivalence and subject to some mild ``positivity" conditions) on $X$ and a positive integer $d$
we have constructed the scheme $\JA$ which fits the above framework. Set-theoretically $\JA$ can be described as the set  of pairs
$(\SEE)$, where $\EE$ is a torsion-free sheaf on $X$ having rank $2$ and the Chern invariants $c_1 (\EE) =L$ and $c_2 (\EE) =d$,
and $[e]$ is a projectivized\footnote{i.e. the homothety equivalence class of a global section $e$ in $H^0 (\EE)$.} global section of $\EE$ whose zero-locus $Z_e = (e=0)$ is a subscheme of dimension $0$ in $X$. From this description it
follows that $\JA$ admits natural maps to the Hilbert scheme $X^{[d]}$, the scheme parametrizing the subchemes of $X$
having dimension $0$
 and length $d$, and the moduli stack $M_X (2, L,d)$ of torsion free sheaves on $X$ having rank $2$ and the Chern invariants
$L$ and $d$.
 This is encapsulated in the following diagram
\begin{equation}\label{HJM}
\xymatrix{
&{\JA} \ar[dl]_h \ar[dr]^{\pi}& \\
  {{\bf M}_X (2,L,d)  }& &{\XD} }
\end{equation}
where the morphism $h$ (resp. $\pi$)  takes a pair $(\SEE)$ to $[\EE]$, the point of the moduli stack  $M_{X} (2,L,d)$
 corresponding to the sheaf $\EE$ (resp. $[Z_e]$, the point of the Hilbert scheme $\XD$ corresponding to the subscheme
$Z_e = (e=0)$).
\\
\indent
Thus $\JA$ can be thought of as a kind of thickening of  $\XD$ and $M_{X} (2,L,d)$. In both cases the thickening is obtained by
inserting rational fibres. The fibre $U_{\EE}$ of $h$ over a point $[\EE] \in M_{X} (2,L,d)$ is easily seen to be the Zariski open
subset of $\PP(H^0 (\EE))$\footnote{ For a nonzero vector space $V$ we denote by $\PP(V)$ the projective space associated to $V$, which is
the space of one-dimensional subspaces of $V$.} parametrizing global sections of $\EE$ whose zero-locus is zero-dimensional. 
The fibre of $\pi$ over $[Z] \in \XD$ will be denoted $\JAA_{Z}$. This is the set of all torsion-free sheaves of rank $2$ on $X$
with Chern invariants $L,d$ and a global section whose zero-locus is the subscheme $Z$. To see what this set looks like, observe that
if $\EE$ is a sheaf in $\JAA_{Z}$ and if $e$ is a global section of $\EE$ with $(e=0)=Z$, the sheaf $\EE$ fits into the following
exact sequence of sheaves, the Koszul sequence defined by the global section $e$ of $\EE$ 
\begin{equation}\label{exs}
\xymatrix@1{
0\ar[r]& {\OO_X} \ar[r]& {\EE}\ar[r]&{\ID_Z} (L) \ar[r] & 0 }
\end{equation}
where $\ID_Z$ is the sheaf of ideals of $Z$ in $X$, the monomorphism is defined by the multiplication by $e$ and the epimorphism is
defined by taking the exterior product with $e$.
\\
\indent 
This gives an identification of $(\EE,e)$ with an element $\alpha$ of the group of extensions \linebreak $Ext^1 (\ID_Z (L), \OO_X)$.
Multiplying $e$ by a non-zero scalar will multiply $\alpha$ by the same scalar. Thus $\JAA_{Z}$ is the projective space
$\PP(Ext^1 (\ID_Z (L), \OO_X))$.
\\
\indent
Analyzing $\JA$ as a fibre space over $\XD$ we show in \cite{[R1]} that $\JA$ provides a very natural passage from
configuration of points on $X$ to realm of
\\
\indent
{\it 1) reductive Lie algebras and their representations
\\
2) Fano toric varieties, whose hyperplane sections are Calabi-Yau varieties}.
\\
\\
One of the consequences of  1) is that the methods of the representation theory of Lie algebras can be applied to study 
algebro-geometric properties
of configurations of points on smooth complex projective surfaces, while 2) should provide new invariants of smooth
complex projective surfaces based on invariants of Fano toric and Calabi - Yau varieties.
\\
\\
\indent
One of the key properties of $\JA$ leading to 1) and 2) above is a surprisingly rich algebraic structure hiding  behind
classical configurations of points on surfaces. More precisely, for suitable points $[Z] \in \XD$, where $Z$ is a subscheme
of $d$ distinct points  on $X$ (such subschemes of $X$ in the sequel will be called {\it configurations} of points on $X$),
the space of complex valued functions $H^0 (\OO_Z) $ on $Z$ carries a sort of Hodge structure.
This will be expounded upon in more details in {\bf Part II} of this introduction but let us pose for a moment and explain what we mean by that.
\\
\indent
It is well-known that for a compact complex manifold $Y$ equipped with a K\"ahler form the DeRham cohomology groups
$H^m (Y,{\CC}),\,\,m=0,\ldots,2dim_{\CC} Y$, admit the direct sum decomposition
\begin{equation}\label{HD}
H^m (Y,{\CC}) = \bigoplus_{p+q =m} H^{p,q} (Y)
\end{equation}
into the Hodge groups $H^{p,q} (Y)$.
\\
\indent
In our case the role of $Y$ is taken over by configurations $Z$ of $d$ distinct points on $X$ and the only nontrivial De Rham 
cohomology group
$H^0 (Z, {\CC}) = H^0 (\OO_Z) $ has no nontrivial Hodge decomposition in the classical sense. But it turns out that equipping $Z$ with
a suitable extension class $[\AL]$ in $\JAA_{Z} = \pi^{-1} ([Z])$, the fibre of $\pi$ in (\ref{HJM}), gives rise to a distinguished direct sum 
decomposition
\begin{equation}\label{HD1}
H^0 (Z,{\CC}) = H^0 (\OO_Z) = \bigoplus^l_{p=0} H^{p} (Z,[\AL])
\end{equation}
where the summands $H^{p} (Z,[\AL])$ as well as the weight $l$ of this decomposition have precise geometric meaning (see Remark \ref{geomm}).
\\
\indent 
By analogy with the classical theory one can think of $[\AL] \in \JAA_{Z}$ as a non-classical K\"ahler structure of $Z$ which induces
the above Hodge-type decomposition of the De Rham group $H^0 (Z,{\CC})$ of $Z$. Furthermore, as $[\AL] $ varies in some non-empty Zariski open
subset of $\JAA_{Z}$ one obtains the variation of the decomposition in (\ref{HD1}). This variation, as the classical Variation of Hodge
structure, possesses Griffiths transversality property (see \cite{[G]} and \S4.3, \cite{[R2]}, for more details).
Putting together the decomposition (\ref{HD1}), its Griffiths transversality property and the multiplicative structure of $H^0 (\OO_Z)$
one arrives to a reductive Lie algebra intrinsically associated to $\ZA$.
\\
\indent 
Heuristically, the apparence of this Lie algebra means  that a classical static configuration of points $Z$ becomes a dynamical object and by going from
$Z$ to the fibre $\JAA_Z$ one sees ``quantum" phenomena. For example, the classical observables, functions 
on $Z$, can be turned into linear operators on $H^0 (\OO_Z) $.
 In particular, to every point $z \in Z$ one can associate such linear operators and speak about 
interactions between distinct points of $Z$. Furthermore, the above mentioned Hodge structure and linear operators can be
conveniently encoded in a certain trivalent graph and one can define a sort of ``vibration" of the linear operators along this 
trivalent graph. This ``vibration", interpreted geometrically, yields Calabi-Yau varieties associated to points of $Z$. Thus 
 from the rigorous mathematical constructions emerges a picture evocative of string theory.
\\
\\
\indent
Hopefully this rather sketchy recall of \cite{[R1]} gives an idea that the non-abelian Jacobian $\JA$ is an interesting object 
establishing natural bridges between Hilbert schemes of points and vector bundles on algebraic surfaces, reductive Lie algebras, Calabi-Yau varieties.
In addition, all the constructions are permeated by very suggestive relations with string theory. Given all this we thought that it might be useful
to supplement \cite{[R1]} with the present paper whose goals are:
\\
a) explain the main mathematical aspects of 1) and 2) above for a variety of any dimension $\geq 2$
\\
b) draw analogies of certain aspects of our mathematical constructions with string theory.
\\
\\
For b), a reader, who is an expert in string theory, should be cautioned: the author is a total neophyte in physics,
so our analogies might strike such a reader as very naive and superficial. However, these analogies are so suggestive that
we thought it might be of some interest to physicists.
\\
\\
\indent
In the rest of the introduction we will outline the main results and ideas of the paper for configurations coming from holomorphic vector bundles or locally free sheaves (as it is customary we make no distinction between these notions). This is largely for two reasons. The first one is that
these configurations possess technical properties which are necessary for our constructions.
Once these technical properties are isolated a passage to a more general setting presents no
difficulties. The second reason is that the vector bundles constitute one of the obvious
targets for applications of our constructions.
\begin{center}
{\bf
Outline of the results in a vector bundle setting}
\end{center}
Let $X$ be a smooth complex projective variety of dimension $n \geq 2$
(all dimensions are over $\CC$, unless stated otherwise) and let
$\EE$ be a holomorphic vector bundle of the same rank $n$ over $X$. Denote by
$\OL$ the determinant bundle of $\EE$ and assume
$H^i (\OO_X (-L)) = 0$, for all $i\leq n-1$, e.g.
$\OL$ is an ample line bundle. This is somewhat restrictive assumption but it holds 
in many geometric situations which motivated our considerations.
\\
\indent
In addition, it always will be assumed that $\EE$ has a global section $e$ whose zero-locus
$Z_e =(e=0)$ is a set of distinct points, i.e. the section $e$ has finite number of zeros and every zero of $e$ has multiplicity 1. These ``good" sections of $\EE$ are considered equivalent if they differ by a nonzero scalar multiple. Their equivalence classes form a 
Zariski open subset $\RU_{\EE}$ of the projective space $\PP(H^0(\EE))$. Over every point
$[e] \in \RU_{\EE}$ we have its zero-locus $Z_e =(e=0)$ and as $[e]$ varies in
$\RU_{\EE}$ the zero-loci $Z_e$ will move in $X$ and form the family
\begin{equation}\label{pue}
p_{\RUE} : \ZDO \longrightarrow \RU_{\EE} 
\end{equation} 
of 0-dimensional subvarieties of $X$. Observe that $p_{\RUE}$ is an unramified covering whose degree is the $n$-th Chern number of $\EE$, i.e. 
\begin{equation}\label{d}
d: = deg(p_{\RUE} )= \int_{[X]} c_n (\EE)
\end{equation} 
is the integrated $n$-th Chern class $c_n (\EE)$ of $\EE$.
\\
\indent
 The variety $\ZDO$ is a natural geometric object
associated with $(X,\EE)$ and it will be called the 
{\it universal configuration} of $\EE$. The open set $\RU_{\EE}$ is the parameter space for these configurations. We suggest that this is only ``visible" or ``classical"
parameters of the point configurations and $\ZDO$ itself is only a ``visible" part of the
``universe" associated to $\EE$. A construction of the space of ``hidden" parameters for
point configurations $Z_e, [e] \in \RU_{\EE}$, and uncovering its various features constitute the mathematical content of the paper.
\begin{center}
{\bf 
Part I: The Jacobian of $\RUE$ or the space of ``hidden" parameters}
\end{center}
The main result of this part is a construction of the torsion free sheaf
$\EXTSO$ over $\RUE$ whose fibre $\EXTSO ([e])$ at $[e] \in \RUE$ is isomorphic to the group of extensions $Ext^{n-1} (\ID_{Z_e} (L), \OO_X)$, where $\ID_{Z_e}$ is the sheaf of ideals of 
$Z_e$. In the sequel this group will be denoted by $\EXT$.
\\
\\
\indent
Let $\delta(\RUE)$ be the rank of $\EXTSO$ and let $Sing(\EXTSO)$ be the singularity set of
$\EXTSO$, i.e. the subset of $\RUE$ where the rank of fibres of $\EXTSO$ exceeds 
$\delta(\RUE)$. From general properties of torsion free sheaves (see e.g. \cite{[O-S-S]})
it follows that $Sing(\EXTSO)$ is a closed proper subscheme of $\RUE$.
Set $\RUEP = \RUE \setminus Sing(\EXTSO)$. The restriction of 
$\EXTSO$ to $\RUEP$ is now a locally free sheaf of rank $\delta(\RUE)$. In particular,
$dim(\EXT) = \delta(\RUE)$, for every $[e] \in \RUEP$.
\\
\indent
Define 
\begin{equation}\label{jac}
\JUEO =\PP(\EXTSO \otimes \OO_{\RUEP})
\end{equation}
 to be the projectivized bundle associated
to $\EXTSO \otimes \OO_{\RUEP}$. This is our space of ``hidden" parameters for the configuration $\ZDO$ or what we call {\it the nonabelian Jacobian} of $\RUE$.
\\
\indent
By definition our Jacobian comes with the projection
\begin{equation}\label{pi}
\pi_{\RUE} : \JUEO \longrightarrow \RUEP
\end{equation}
and the line bundle $\OO_{\JUEO} (1)$  defined so that its direct image
$$
\pi_{\RUE \ast} \left(\OO_{\JUEO} (1) \right) = \left(\EXTSO \otimes \OO_{\RUEP} \right)^{\ast}
$$
Thus the construction of the space of ``hidden" parameters amounts to attaching the projective space $\PP(\EXT)$ to every point $[e] \in \RUEP$ of the space of ``visible" parameters of the 
universal configuration of $\EE$. It is clear that we get something new only if the rank
$\delta(\RUE) \geq 2$. This always will be assumed unless said otherwise.
Observe that for $n=1$ our construction {\it never} gives anything new. So an existence of nontrivial
nonabelian Jacobian is a distinguished feature of higher dimensional geometry.
\\
\indent
As the space of ``visible" parameters $\RUEP$ governs the ``classical" motion of our
configurations in $X$, we expect that the ``hidden" parameters in $\JUEO$ will reveal
some ``quantum" properties of the universal configuration $\ZDO$. This is indeed the case
and the feature of $\JUEO$ responsible for this is a Hodge-like structure and its variation
along the fibres of $\JUEO$.
\begin{center}
\bf
Part II: A variation of Hodge-like structure associated to $\JUEO$
\end{center}
We begin by explaining what happens on the fibre of $\JUEO$ over a point $[e] \in \RUEP$.
On the ``classical" level over $[e]$ lies the configuration $Z_e$ of $d$ distinct points
(particles) on $X$. The space of functions $\HOZ$ on $Z_e$ can be viewed as the space of
classical observables of $Z_e$. Attaching to $[e]$ the ``hidden" parameters
$\PP(\EXT)$ has the following effect.
\\
\indent
There is a nonempty Zariski open subset $\JBB([e])$ of $\PP(\EXT)$ such that for every
$[\alpha]\in \JBB([e])$ the space of classical observables $\HOZ$ admits a distinguished
direct sum decomposition
\begin{equation}\label{dsd}
\HOZ = \bigoplus^{w(\RUE)}_{p=0} \HB^p \EA
\end{equation}
where the number of summands $w=w(\RUE) \geq 1$ is called the weight of $\RUE$.
\\
\indent
The summand $\HB^0 \EA$ plays a special role.
\begin{lem}\label{sum0}
\begin{enumerate}
\item[1)]
$\HB^0 \EA$ contains the subspace $H^0(\OO_X) =\CC$ of constant functions on $Z_e$.
\item[2)]
There is a natural isomorphism of vector spaces
\begin{equation}\label{iso}
\theta_{[e]} (\alpha) : \HB^0 \EA \longrightarrow \EXT
\end{equation}
which sends the constant function $1 \in \HB^0 \EA$ to the extension class
$\alpha \in \EXT$.
\end{enumerate}
\end{lem}
The essential point is that the decomposition (\ref{dsd}) varies with $\EA$.
In particular, we study its variation with respect to the ``hidden" parameter
$[\alpha]$. By virtue of the identification in (\ref{iso}) this variation is closely related to the multiplication in the ring $\HOZ$ by elements in $\HB^0 \EA$. It turns out that the
multiplication by $t \in \HB^0 \EA$ preserves the subspace
\begin{equation}\label{HT-w}
\HT_{-w} \EA = \bigoplus^{w-1}_{p=0} \HB^p \EA
\end{equation}
and acts on this decomposition by shifting the grading (up or down) at most by 1.
In other words, if $D(t)$ denotes the operator of the multiplication by $t$ in
$\HOZ$ and if $D_p (t)$ is its restriction to the summand $\HB^p \EA$, for
$p=0,\ldots,w-1$, then
\begin{equation}\label{Dpt}
D_p (t) : \HB^p \EA\longrightarrow \HB^{p-1} \EA \oplus \HB^p \EA \oplus \HB^{p+1} \EA
\end{equation}
with an understanding that $\HB^{-1} = 0$.
Denote by $D^{-}_p (t) $ (resp. $D^{0}_p (t) $ and $D^{+}_p (t) $) the projection of
$D_p(t)$ on the first (resp. the second and the third) summand in (\ref{Dpt}) with an understanding that $D^{-}_0 (t) = D^{+}_{w-1} (t) = 0$. 
\\
\indent
Set
\begin{equation}\label{D+0-}
D^{\pm} (t) = \sum^{w-1}_{p=0} D^{\pm}_p (t),\,\,\,\,D^{0} (t) = \sum^{w-1}_{p=0} D^{0}_p (t)
\end{equation}
The operator $D(t)$ restricted to the subspace $\HT_{-w} \EA$ in (\ref{HT-w}) admits the following decomposition
\begin{equation}\label{Ddec}
D(t) = D^{-} (t) + D^{0} (t) +D^{+} (t)
\end{equation}
where $D^{0} (t)$ preserves the grading in (\ref{HT-w}) and 
$D^{\pm} (t)$ shifts it by $\pm1$. Thus the ``hidden" variable $[\alpha]$ produces the direct sum
decomposition (\ref{dsd}) and defines the subspace
\begin{equation}\label{esp}
\{ D^{\pm} (t),D^{0} (t) \mid t\in \HB^0 \EA \}
\end{equation}
of $End(\HT_{-w} \EA)$. Defining the Lie algebra 
\begin{equation}\label{Lie}
{\bf \goth{g}}\EA
\end{equation}
  as the Lie subalgebra of $End(\HT_{-w} \EA)$ generated by the subspace in (\ref{esp}) we obtain a Lie algebra intrinsically associated
to every point $\EA$ in an appropriate non-empty Zariski open subset of $\JUEO$.
\\
\indent 
Heuristically, the above considerations imply that {\it a priori} static finite set of points is turned into
a dynamic object. The dynamics manifests itself in a variation of the decomposition
(\ref{HT-w}) with respect to the ``hidden" variable $[\alpha]$. The variation itself is governed by the decomposition $(\ref{Ddec})$. In particular, we see that our classical observable,
a function $t \in \HB^0 \EA$, acquires now the status of a system of linear operators and while
the multiplication operators $D(t)$ are obviously commuting the Lie subalgebra of
$End(\HT_{-w} \EA)$ generated by the operators in (\ref{esp}), in general, is not abelian.
All this indicates that our ``hidden" parameter $[\alpha]$ is capable of detecting
``quantum" properties of the configuration $Z_e$.
\\
\indent
Until now we have been discussing the situation on the fibre of $\JUEO$ over a point
$[e]$ in $\RUEP$. What happens when we vary $[e]$? First of all the rings
$\HOZ$ fit together to form the locally free sheaf $p_{\RUE \ast} \OO_{\ZDO}$, the direct image
of $\OO_{\ZDO}$ under the morphism $p_{\RUE}$ in (\ref{pue}). So to generalize our results we
need to consider the pullback
$\PIPO$ of $p_{{\RUE} \ast} \OO_{\ZDO}$ via the projection $\pi_{\RUE}$ in (\ref{pi}). We can now state the sheaf version of our discussion.
\begin{thm}\label{Th1}
\begin{enumerate}
\item[1)]
There is a distinguished nonempty Zariski open subset $\JBBU$ of $\JUEO$ and an integer
$w(\RUE) \geq 1$ such that the sheaf
$\PIPO \otimes \OO_{\JBBU}$ admits a distinguished direct sum decomposition
\begin{equation}\label{dsds}
\PIPO \otimes \OO_{\JBBU} = \bigoplus^{w(\RUE)}_{p=0} \HB^p_{\RUE}
\end{equation}
The integer $w=w(\RUE)$ is called the weight of $\RUE$ (if no ambiguity is likely we will 
simplify the above notation by omitting the reference to $\RUE$).
\item[2)]
The summand $\HB^0$ in (\ref{dsds}) has the following properties.
\begin{enumerate}
\item
$\OO_{\JBB}$ is a subbundle of $\HB^0$.
\item
There is a natural isomorphism
\begin{equation}\label{theta-s}
\theta: \HB^0 \otimes \OO_{\JBB} (-1) \longrightarrow \pi^{\ast}_{\RUE} \EXTSO
\end{equation}
which induces the following commutative diagram
\begin{equation}\label{theta-d}
\xymatrix{
0 \ar[d]& 0\ar[d]\\
\OO_{\JBB} (-1) \ar[d] \ar@{=}[r]& \OO_{\JBB} (-1)\ar[d] \\
\HB^0 \otimes \OO_{\JBB} (-1) \ar[d] \ar[r]^{\theta} &\pi^{\ast}_{\RUE} \EXTSO \ar[d]\\
\left(\HB^0 / \OO_{\JBB} \right) \otimes \OO_{\JBB} (-1) \ar[d] \ar[r]^(.6){\theta^{\prime}} &
\RELTAN (-1) \ar[d] \\
0&0 }
\end{equation}
where $\RELTAN$ is the relative tangent bundle of
$\pi_{\RUE}: \JBB \longrightarrow \RUE$ and the column on the right is the relative Euler sequence on $\JBB$ tensored with $\OO_{\JBB} (-1)$.
\end{enumerate}
\end{enumerate}
\end{thm}
Next we state the sheaf version of the endomorphisms in (\ref{esp}).
\begin{thm}\label{Th2}
Set
\begin{equation}\label{HT-wS}
\HT_{-w} = \bigoplus^{w-1}_{p=0} \HB^p
\end{equation}
\begin{enumerate}
\item[1)]
The multiplication in the sheaf of rings $\PIPO$ gives rise to the following morphism of sheaves
\begin{equation}\label{D}
D: \HT_{-w} \longrightarrow \left(\HB^0\right)^{\ast} \otimes \HT_{-w}
\end{equation}
which admits the decomposition
\begin{equation}\label{Ddec-m}
D= D^{-} + D^{0} + D^{+}
\end{equation}
where $D^0,D^{\pm} \in Hom (\HT_{-w}, \left(\HB^0\right)^{\ast} \otimes \HT_{-w})$
\item[2)]
The morphism $D^0$ (resp. $D^{\pm}$) is a morphism of degree $0$ (resp. $\pm1$) with respect to
the grading in (\ref{HT-wS}). Furthermore the morphisms
$D^0,D^{\pm}$ are subject to the following relations
\begin{equation*}
(i)\,\,\, D^{\pm} \wedge D^{\pm} =0 
\end{equation*}
\begin{equation}\label{rel}
(ii)\,\,\,D^0 \wedge D^{\pm} + D^{\pm} \wedge D^0 =0 
\end{equation}
\begin{equation*}
(iii)\,\,\,D^0 \wedge D^0 + D^{+}\wedge D^{-} + D^{-} \wedge D^{+}=0 
\end{equation*}
\item[3)]
The morphisms $D^0,\,D^{\pm}$ define the sheaf ${\bf\goth{G}}_{\JBBU}$ of reductive Lie algebras on $\JBBU$. The fibre
of ${\bf\goth{G}}_{\JBBU}$ at every point $\EA \in \JBBU$ is the Lie algebra ${\bf\goth{g}} \EA$ defined in (\ref{Lie}).
\end{enumerate}
\end{thm}
\begin{center}
{\bf Part III: A geometric interpretation of Theorem \ref{Th2}  }
\end{center}
The commutativity of the multiplication in $\PIPO$ translates into the fact
that $D^2= D \wedge D=0$, i.e. $D$ is a Higgs morphism of $\HT_{-w}$
with values in $\left(\HB^0 \right)^{\ast}$ (see \cite{[S]} or Definition 3.2,\cite{[R1]},
for definitions). To extract geometry from Theorem \ref{Th2} we deform the morphism $D$ using the direct sum decomposition (\ref{HT-wS}). This gives us a particular family of Higgs morphisms
parametrized by an affine variety explicitly defined by a system of quadratic equations 
(Lemma \ref{H-hat}). This variety is denoted by ${\hat{H}}_{\RUE}$. Passing to the homothety equivalence classes of
Higgs morphisms parametrized by ${\hat{H}}_{\RUE}$ we obtain the quotient
\begin{equation}\label{alb}
H_{\JBBU} =  {\hat{H}}_{\RUE} / {\bf \CC^{\ast}}
\end{equation}
which is called the nonabelian Albanese of $\JBBU$.
\begin{thm}\label{Th3}
Let the weight $w=w(\RUE) \geq 2$. Then $H_{\JBBU}$ is a projective toric variety of dimension
$w-1$. The variety comes equipped with a very ample invertible sheaf 
$\OO_{ H_{\JBBU}} (1)$ such that the dualizing sheaf 
$\omega_{ H_{\JBBU}} = \OO_{ H_{\JBBU}} (-1)$, i.e. the variety
$H_{\JBBU}$ is a (toric) Fano variety. In particular, if the weight
$w\geq 3$, then the hyperplane sections of $H_{\JBBU}$ are Calabi-Yau varieties of dimension
$w-2$.
\end{thm}
The nonabelian Albanese of $\JBBU$ depends only on the weight $w$ and the relations between the morphisms $D^0,D^{\pm}$ in (\ref{rel}). So it completely ``forgets" about the variety $X$
and the universal configuration $\ZDO$. However, there are various ways to remedy this by constructing correspondences between $\JBBU$ and $H_{\JBBU}$. In particular, there is a geometric correspondence described as follows.
\begin{thm}\label{Th4}
There is a natural geometric correspondence between $\JBBU$ and $H_{\JBBU}$ given by the
morphism
\begin{equation}\label{cy}
CY_{\RUE}: \JBBU \longrightarrow \PP(H^0 (\OO_{ H_{\JBBU}} (d)))
\end{equation}
where $d$ is the $n$-th Chern number of $\EE$. This morphism sends a point
$\EA$ of $\JBBU$ to a cycle of the hyperplane sections of $H_{\JBBU}$. More precisely,
the (projectivized) section $ CY_{\RUE} \EA$ corresponds to the divisor
\begin{equation}\label{Hz}
\sum_{z\in Z_e} H_{z,[\AL]}
\end{equation}
where $H_{z,[\AL]}$ is a hyperplane section of $H_{\JBBU}$ intrinsically associated
to a point $(z,[\AL]) \in (Z_e,[\AL])$.
\end{thm}
 If $w \geq 3$, then every $H_{z,[\AL]}$ is a 
Calabi-Yau variety and the morphism $CY_{\RUE}$ associates to every point of $\JBBU$ a cycle
of Calabi-Yau varieties. For this reason $CY_{\RUE}$ is called Calabi-Yau cycle map.
\\
\\
\indent
The result of Theorem \ref{Th4} tells us that once we fix $[\AL]$ in the fibre of $\JBBU$ over a point $[e]$ in $\RUEP$, every point $z\in Z_e$ ``opens up" and reveals the Calabi-Yau
variety $H_{z,[\AL]}$. Thus our constructions show that over every point of the ``visible
universe" $\ZDO$ ``seats" a Calabi-Yau variety of dimension $w-2$. These Calabi-Yau varieties
have moduli: for every ``classical" configuration $Z_e$ this moduli is our 
``hidden" parameter $[\AL]$ moving in the open subset $\JBBU \cap \PP(\EXT)$.
Furthermore the hyperplane sections $H_{z,[\AL]}$ in (\ref{Hz}) arise naturally as a kind of 
propagation of the $\delta$-function $\delta_z$ supported at $z\in Z_e$. This 
propagation is provided by the morphisms 
$D^0,D^{\pm}$ in Theorem \ref{Th2}. The main point here is that the data of the decomposition
(\ref{HT-wS}) and the morphisms $D^0,D^{\pm}$ can be organized in a trivalent graph which is
denoted $G(\JBBU)$ (see (\ref{G})). So the above mentioned propagation is obtained by following the paths of the graph $G(\JBBU)$ with endomorphisms from the space (\ref{esp}). This way one obtains
``strings"=strings of operators intrinsically attached to every point
$(z,[\AL]) \in (Z_e,[\AL])$. In particular, the creation of 
$H_{z,[\AL]}$ could be regarded as a ``vibration" of a string of operators according to the pattern of the
graph $G(\JBBU)$. 
\\
\indent
Thus what emerges from our results is a picture of ``universe" = ``visible"+``hidden"
associated to the universal configuration $\ZDO$ which is very suggestive of the one
predicted by string theory. From this perspective our results could be succinctly formulated
as follows.
\\
\\
{\bf Meta-principle.}
{\it All smooth complex projective varieties of dimension $\geq 2$ have
a string theoretic nature:
\\
$\bullet$ points on these varieties have ``hidden" parameters
\\
$\bullet$ in the
space of these ``hidden" parameters points are no longer 0-dimensional but
become ``strings" in a precise mathematical sense
\\
$\bullet$ these ``strings" ``vibrate" according to the pattern of a certain trivalent
graph and this vibration creates Calabi-Yau varieties.}
\\
\\
\indent
In the main body of the paper we tried to convey principal ideas of our constructions with a bear minimum of technicalities.
For more detailed discussion and complete proofs an interested reader is
 referred to \cite{[R1]}. This reference contains most of the ideas
exposed here. The only really new result of this paper is a generalization of the results in
\cite{[R1]} to an arbitrary
dimension $\geq 2$.
\\
\indent
The rest of the paper is organized as follows.
\\
- In \S1 definitions of the universal cluster and the universal configuration are given.
\\
- In \S2 we define the space of ``hidden" parameters of the universal cluster.
\\
- In \S3 the Hodge-like structure alluded to in {\bf Part II} is defined.
\\
- In \S4 we discuss mathematics of theorems stated in {\bf Part II} 
\\
-\S\S4.1,4.2 contain proofs of Theorem \ref{Th1} and Theorem \ref{Th2} respectively,
\\
-\S4.3 treats a trivalent graph which is an explicit mechanism for transforming
classical observables=functions to ``quantum" observables=path-operators,
\\
-{\bf Part III} is treated in \S5.
\\
-\S6 is entirely devoted to speculations about possible string theoretic nature of our constructions.
\\
- In \S7 we spell out conditions on a family of configurations of points of $X$ for which all our constructions
remain valid. This allows to generalize all our results to families of configurations which 
{\it do not}
necessarily come from vector bundles. This becomes important in dimension $\geq 3$,
where one no longer has an equivalence between extension classes and torsion free sheaves
(as it is the case in dimension 2).
\\
-In \S8 we illustrate our general considerations by an example of configurations which are complete intersections.
\\
\\
\indent
It is a pleasure to thank Vladimir Roubtsov for his interest and encouragement.

\section{The universal configuration or the ``visible universe" associated to $\EE$}

Let $\EE$ be a holomorphic vector bundle of rank $n$ over an $n$-dimensional ($n\geq2$)
smooth complex projective variety $X$. The basic topological invariants of
$\EE$ are its Chern cohomology classes 
$c_i (\EE) \in H^{2i} (X,{\bf Z}),i=0,\ldots,n$. Under some favorable circumstances the
Poincar\'e duals of $c_i (\EE)$ can be represented by algebraic subvarieties or, more
generally, by subschemes of $X$. We will be particularly interested in the top Chern class
$c_n (\EE)$. It is well-known (see e.g. \cite{[G-H]}) that if $\EE$ has a global section,
say $e$, such that its zero-locus $Z_e =(e=0)$ is 0-dimensional, then $Z_e$ is a representative of the Poincar\'e dual of $c_n (\EE)$. In particular, the degree
$d=deg(Z_e)$, i.e. the number of zeros of $e$ counted with their multiplicities, is given by
\begin{equation}\label{deg}
d=deg (Z_e) =\int_{[X]} c_n (\EE)
\end{equation}
the Gauss-Bonnet type formula.
\\
\indent
One can view $Z_e$ as a geometric realization of $c_n (\EE)$. From now on we assume that the 
vector bundle $\EE$ is subject to the following condition
\begin{equation}\label{s-cond}
H^0(\EE) \neq 0\,\,\,and\,\,\,for\,\,\,some\,\,\,e\in H^0(\EE)\,\,\,the\,\,\,
zero-locus\,\,\,Z_e =(e=0)\,\,\, is \,\,\,0-dimensional
\end{equation}
This implies in particular that $d=deg(Z_e) >0$.
\\
\indent
The projective space $\PP(H^0(\EE))$ parametrizes the zero-loci $Z_e$ as $e$ moves
in $H^0(\EE) \setminus \{0\}$. We distinguish in $\PP(H^0(\EE))$ the subset 
$\UE$ corresponding to the sections of $\EE$ with 0-dimensional zero-loci.
This is a Zariski open subset of $\PP(H^0(\EE))$. The assumption (\ref{s-cond}) implies that
$\UE$ is nonempty.
\\
\indent
Consider the incidence correspondence
\begin{equation}\label{inc}
X\times \UE \supset \ZD = \{(x,[e]) \in X\times \UE \mid e(x)=0\}
\end{equation}
From the definition of $\UE$ it follows that the projection
\begin{equation}\label{p2}
p_2 : \ZD \longrightarrow \UE
\end{equation}
on the second factor of $X\times \UE$ is a finite morphism of degree $d$. Thus 
$\ZD$ can be characterized as the largest family of 0-dimensional subschemes of $X$ which
provide geometric realization of $c_n(\EE)$. We will call it 
{\it the universal cluster}\footnote{A $0$-dimensional subscheme of $X$ is also called {\it cluster} of $X$. This is
the terminology proposed by Miles Reid about 20 years ago.} \label{clust} of $\EE$ or the geometric $n$-th Chern class of $\EE$.
\\
\indent
The branch locus of the finite morphism $p_2$ in (\ref{p2}) corresponds to sections
of $\EE$ with some of its zeros being multiple.
\begin{defi}\label{s-reg}
A global section $e$ of $\EE$ is called regular if its zero-locus $Z_e$ is a finite set
and every zero of $e$ has multiplicity $1$.
\end{defi}
From now on we assume that $\EE$ has a regular section and we denote by
$\RUE$ the subset of $\UE$ corresponding to the regular sections of $\EE$, i.e.
$\RUE$ is the set of regular  ``values" of $p_2$. Denote by $\ZDO$ the inverse
image of $\RUE$ with respect to the morphism $p_2$ in (\ref{p2}). Then its restriction
\begin{equation}\label{pUE}
p_{\RUE} : \ZDO \longrightarrow \RUE
\end{equation}
is an unramified covering of degree $d$.
\begin{defi}\label{conf}
The variety $\ZDO$ is called the universal configuration associated to $\EE$
(or simply the universal configuration of $\EE$).
\end{defi}
\begin{rem} The fact that $\RUE$ is nonempty is certainly a restrictive assumption.
It holds, for example, for vector bundles generated by global sections.
Such bundles arise naturally in geometric situations, e.g. from a (non-constant)
morphism of $X$ into the Grassmannian $Gr(n,V)$ of $n$-planes of a given complex
vector space $V$. So vector bundles with regular sections should be viewed as
``geometric" bundles.
\end{rem}
The variety $\ZDO$ could be thought of as a configuration of points (particles)
moving in $X$. The open set $\RUE$ is the space of ``visible" or ``classical"
parameters of this configuration. We now turn to a construction of the space of
its ``hidden" parameters.

\section{The sheaf of extensions $\EXTS$ and the space of ``hidden" parameters
associated to $\ZD$}
In this section the construction outlined in {\bf Part I} of the introduction will be given.
To avoid unnecessary technical complications it will be assumed that the determinant bundle
$\OL$ of $\EE$ is subject to the following
\begin{equation}\label{L-cond}
H^i (\OO_X (-L)) = 0,\,\,\forall i\leq n-1
\end{equation}
This condition holds, for example, if $\OL$ is ample.
\\
\indent
Our goal is to construct a certain sheaf $\EXTSO$ over $\RUE$ whose fibre $\EXTSO ([e])$ over a point $[e] \in \RUE$
is isomorphic to the group $Ext^{n-1} (\ID_{Z_e} (L), \OO_X)$.
\\
\indent
There are some simple algebro-geometric reasons for considering the extension groups
\begin{equation}\label{Ext}
\EXT = Ext^{n-1} (\ID_{Z_e} (L), \OO_X)
\end{equation}
where $[e]\in \UE$ and $\ID_{Z_e}$ is the ideal sheaf of $Z_e$. Let us  discuss them briefly.
\\
\indent
It is well-known fact that $Z_e$ fails to impose independent conditions on the linear system
$\LKS$, where $K_X$ is the canonical divisor of $X$.
For points to impose (or not) independent conditions on a linear system is a classical
algebro-geometric notion. To put it in modern terms in our context consider the short
exact sequence of sheaves
\begin{equation}\label{s1}
\xymatrix{
0\ar[r]&\ID_{Z} (\LK)\ar[r]& \OO_X (\LK)\ar[r]^(.5){\rho_Z}&\OO_{Z} (\LK) \ar[r]&0 }
\end{equation}
where $Z$ is a 0-dimensional subscheme of $X$ and $\ID_{Z}$ is its sheaf of ideals.
  \begin{defi}\label{basic-def}
   \begin{enumerate}
   \item [(i)]
  A cluster $Z$ is said to fail to impose independent conditions on
the linear system $\LKS$ or, equivalently and more briefly, $Z$ is called $L$-special , iff 
  the homomorphism
$$
\rho_Z: H^0(\OO_X (\LK)) \longrightarrow H^0(\OO_Z (\LK))
$$
induced from the long cohomology sequence of (\ref{s1}) fails to be surjective.
  \item[(ii)]
  The number $\delta (L, Z) = deg Z \, - rk \,\rho_Z$, where $rk \,\rho_Z$ stands for
the rank of $\rho_Z$, is called the index of 
  $L$-speciality of $Z$.
  \end{enumerate}
  \end{defi}
If we assume that the line bundle $\OO_X (\LK)$ is generated by global sections and hence defines the morphism
$$
\phi_{\LK} : X \longrightarrow \PP(H^0(\OO_X (\LK))^{\ast})
$$
then the index of speciality tells us  that the image of $Z$ under
$\phi_{\LK}$ spans a projective subspace of dimension
$degZ-\delta (L, Z)-1$. So, for example, to say that $2$ points fail to impose independent conditions on the linear system $\LKS$ or, equivalently, $2$ points are $L$-special, means that
these points are mapped by $\phi_{\LK}$ onto the same point,
\\
$3$ points are $L$-special means that the images of these points are on a line,
\\
$4$ points are $L$-special means that the images of these points are on a plane, etc..
\\
Thus it should be clear that the notion of $L$-speciality is a useful concept for studying special projective properties of configuration of points on $X$ with respect to the so called
{\it adjoint} linear system $\LKS$. If we venture into the physical analogy and think of
our point configuration as particles moving in $X$, then $L$-speciality of a configuration
might be regarded as a linear constraint on its motion.
\\
\indent
 We now take $Z=\ZE$ and examine how the notion of its $L$-speciality is related to the extension group $\EXT$ in (\ref{Ext}). This relation becomes clear if
one unravels the long cohomology sequence of (\ref{s1}) beyond the $H^0$-terms. Namely,
one obtains
$$
coker (\rho_{\ZE}) =H^1 (\ID_{\ZE} (\LK))
$$
where we also used our vanishing hypothesis in (\ref{L-cond}). Applying Serre duality
yields
\begin{equation}\label{H1=Ext}
 H^1 (\ID_{\ZE} (\LK))^{\ast} = Ext^{n-1}(\ID_{\ZE} (\LK), \OO_X (K_X)) = 
Ext^{n-1}(\ID_{\ZE} (L),\OO_X)
\end{equation}
Of course the above argument holds for any subscheme of $X$. A particular feature of
$Z_e$ which insures that it is $L$-special is the fact that $\ZE$ comes with a nonzero extension class in 
$\EXT$. This extension class is defined by the Koszul sequence
\begin{equation}\label{Koz}
\xymatrix{
0\ar[r]& \wedge^n \EE^{\ast}\ar[r]^e \ar@{=} [d]&\cdots \ar[r]^e&\wedge^2 \EE^{\ast}\ar[r]^e&
\EE^{\ast}\ar[r]^e&\ID_{\ZE}\ar[r]&0 \\
& \OO_X (-L)&&&&& }
\end{equation}
where $\EE^{\ast}$  is the dual of $\EE$ and the morphisms are contractions with $e$, a section defining $\ZE$. This sequence
is exact precisely when $Z_e$ is 0-dimensional. So (\ref{Koz}) determines a particular
extension class 
$$
\AL_e \in Ext^{n-1} (\ID_{\ZE}, \OO_X (-L)) = Ext^{n-1} (\ID_{\ZE} (L), \OO_X) =\EXT
$$
which is nonzero because $\EE$ is locally free. 
\\
\indent
Thus what comes out of the above discussion is that to every $[e]$ in $\UE$ we can associate
a nonzero complex vector space $\EXT$. These vector spaces fit together to form a sheaf over 
$\UE$. This is the content of the following.
\begin{pro}\label{EXTS}
There exists a torsion free sheaf 
$\EXTS$ over $\UE$ such that its fibre 
$\EXTS ([e])$ is isomorphic to the extension group $\EXT$ in (\ref{Ext}),
 for every $[e]\in \UE$.
\end{pro}
\begin{pf}
Let $\ID_{\ZD}$ be the ideal sheaf of the universal cluster $\ZD \subset X\times \UE$ defined in (\ref{inc}) and consider the defining sequence of $\ZD$
\begin{equation}\label{defZD}
\xymatrix{
0\ar[r]& \ID_{\ZD} \ar[r] & \OO_{X\times \UE} \ar[r] & \OO_{\ZD} \ar[r] & 0 }
\end{equation}
Denote by $p_i$ the projection of $X\times \UE$ onto the $i$-th factor, $i=1,2$.
Taking the pullback $p^{\ast}_1 \OO_X (\LK)$ and tensoring (\ref{defZD}) with it we obtain
\begin{equation}\label{defZD1}
\xymatrix{
0\ar[r]& \ID_{\ZD} \otimes p^{\ast}_1 \OO_X (\LK) \ar[r] &  p^{\ast}_1 \OO_X (\LK)
 \ar[r] & \OO_{\ZD} \otimes p^{\ast}_1 \OO_X (\LK)\ar[r] & 0 }
\end{equation}
Taking the direct image with respect to $p_2$ yields
\begin{equation}\label{defZD2}
\xymatrix{
\scriptstyle{H^0(\LK) \otimes \OO_{\UE}} \ar[r]&
\scriptstyle{p_{2 \ast} \left(\OO_{\ZD} \otimes p^{\ast}_1 \OO_X (\LK) \right)}\ar[r] &
\scriptstyle{R^1p_{2 \ast} \left(\ID_{\ZD} \otimes p^{\ast}_1 \OO_X (\LK) \right)} \ar[r] &0 }
\end{equation}
Dualizing this sequence gives the following
\begin{equation}\label{defZD3}
\xymatrix{
0\ar[r] &
\scriptstyle{\left(R^1p_{2 \ast} \left(\ID_{\ZD} \otimes p^{\ast}_1 \OO_X (\LK) \right)\right)^{\ast}} \ar[r] &
\scriptstyle{\left(p_{2 \ast} \left(\OO_{\ZD} \otimes p^{\ast}_1 \OO_X (\LK) \right)\right)^{\ast}} \ar[r]&
\scriptstyle{H^0(\LK)^{\ast} \otimes \OO_{\UE}} }
\end{equation}
Define
\begin{equation}\label{defEXTS}
\EXTS := \left(R^1p_{2 \ast} \left(\ID_{\ZD} \otimes p^{\ast}_1 \OO_X (\LK) \right)\right)^{\ast}
\end{equation}
From (\ref{defZD3}) it follows that $\EXTS$ is a nonzero subsheaf of the locally free sheaf \linebreak
$\left(p_{2 \ast} \left(\OO_{\ZD} \otimes p^{\ast}_1 \OO_X (\LK) \right)\right)^{\ast}$.
Hence $\EXTS$ is torsion free. Furthermore by definition the fibre of
$\EXTS$ at $[e] \in \UE$ is 
$H^1 (\ID_{\ZE} (\LK))^{\ast}$ which is isomorphic, as we have already seen in (\ref{H1=Ext}), to the extension
group $\EXT$
\end{pf}
\begin{defi} 
The rank of $\EXTS$ is called the index of speciality (with respect to $L$)
of the universal cluster of $\EE$ over $\UE$ and it is denoted $\delta (\UE)$ (or $\delta (L,\UE)$ ).
\end{defi}
Arguing as in {\bf Part I} of the introduction we can define a distinguished
nonempty Zariski open subset $\UE^{\prime}$ over which the sheaf 
$\EXTS$ is locally free of rank $\delta (\UE)$.
\begin{defi}\label{def-jac}
Define
$$
\JUE:= \PP(\EXTS \otimes \OO_{\UE^{\prime}})
$$
the projectivized bundle corresponding to $\EXTS \otimes \OO_{\UE^{\prime}}$.
The variety $\JUE$ is called the nonabelian Jacobian of the universal cluster of $\EE$
over $\UE$.
\end{defi}
The above terminology stems from considerations for surfaces, i.e. $n=2$. In this
case the variety $\JUE$ can be identified with the parameter space of a distinguished family of torsion free
sheaves of rank $2$ over $X$ and its properties evoke a certain analogy with Jacobians of
smooth projective curves (see \cite{[R1]} and \S4, \cite{[R2]}). For $n\geq 3$ an interpretation of $\JUE$ as the parameter space of a family of torsion free sheaves of 
rank $n$ is no longer valid. However we preserved the terminology for lack of a better one.
It is clear that in dimension $n\geq 3$ it is no longer bundles or torsion free sheaves but
rather their complexes (up to a certain equivalence) which are parametrized by $\JUE$.
Thus the variety $\JUE$ can be regarded as the parameter space of a distinguished family
of elements of the derived category of coherent sheaves on $X$. From this perspective
the emergence of our Jacobian  fits naturally with recent understanding of the
important role played by
this category in geometry of projective varieties (see \cite{[B-O]},\cite{[K]}).
\\
\\
\indent
By definition $\JUE$ comes equipped with the natural projection
\begin{equation}\label{proj-pi}
\pi_{\UE} : \JUE \longrightarrow \UE^{\prime}
\end{equation}
and the line bundle $\OO_{\JUE} (1)$ such that its direct image
$$
\pi_{\UE \ast} \OO_{\JUE} (1) = \left(\EXTS \right)^{\ast}
$$
As it was alluded to in the introduction the variety $\JUE$ should play a role of ``hidden"
parameters of the universal configuration $\ZDO$ and it should reveal its ``quantum"
properties. It turns out that one can arrive to this by considering certain Hodge-like
structure on the sheaf $\PIP$.

\section{Hodge-like structure on $\PIP$}
We begin by explaining our construction on a fibre of the projection
$\pi_{\UE} : \JUE \longrightarrow \UE^{\prime}$.
\\
\indent
Fix $[e]\in \UE^{\prime}$. On the fibre of $\JUE$ over $[e]$ the sheaf
$\PIP$ is just the trivial bundle $\HOZ \otimes \OO_{\PP(\EXT)}$.
The Hodge-like structure we have in mind is a distinguished decreasing filtration
on $\HOZ \otimes \OO_{\PP(\EXT)}$. In order to construct it we will need some preliminary
observations.
\begin{lem}
There is a natural morphism
\begin{equation}\label{theta-e}
\theta_{[e]} :\HOZ \otimes \OO_{\PP(\EXT)} (-1) \longrightarrow \left(\HOZK \right)^{\ast}\otimes \OO_{\PP(\EXT)}
\end{equation}
which is an isomorphism at a general point $[\AL] \in \PP(\EXT)$
\end{lem}
\begin{pf}
By Serre duality
\begin{equation}\label{SD1}
\left(\HOZK \right)^{\ast} = Ext^n (\OO_{\ZE} (\LK), \OO_X (K_X))
\end{equation}
The fact that $\ZE$ is 0-dimensional implies that the sheaves
$\EEF^i (\OO_{\ZE} (\LK), \OO_X (K_X)) = 0$, for all $i\leq n-1$. So the local-to-global
spectral sequence for $Ext^n (\OO_{\ZE} (\LK), \OO_X (K_X))$ yields
\begin{equation}\label{SD2}
Ext^n (\OO_{\ZE} (\LK), \OO_X (K_X))= H^0(\EEF^n (\OO_{\ZE} (\LK), \OO_X (K_X)))
\end{equation}
The sheaf $\EEF^n (\OO_{\ZE} (\LK), \OO_X (K_X))$ can be rewritten as follows
\begin{equation}\label{SD3}
\EEF^n (\OO_{\ZE} (\LK), \OO_X (K_X))= \EEF^n (\OO_{\ZE} , \OO_X (K_X))\otimes \OO_X (-L -K_X)=
\omega_{\ZE} \otimes \OO_X (-L -K_X)
\end{equation}
where $\omega_{\ZE} = \EEF^n (\OO_{\ZE} , \OO_X (K_X))$ is the dualizing sheaf of $\ZE$.
This combined with (\ref{SD2}) and (\ref{SD1}) gives
\begin{equation}\label{SD4}
\left(\HOZK \right)^{\ast} = H^0(\omega_{\ZE} \otimes \OO_X (-L -K_X))
\end{equation}
Observe that 
$H^0(\omega_{\ZE} \otimes \OO_X (-L -K_X))$ is an $\HOZ$-module with the structure
multiplication
\begin{equation}\label{stmult}
\HOZ \otimes H^0(\omega_{\ZE} \otimes \OO_X (-L -K_X)) \ARR 
H^0(\omega_{\ZE} \otimes \OO_X (-L -K_X))
\end{equation}
From (\ref{defZD3}) we know that $\EXT$ is a subspace of $\left(\HOZK \right)^{\ast}$.
This combined with (\ref{SD4}) and (\ref{stmult}) yields the pairing
\begin{equation}\label{Ye}
Y_e: \HOZ \otimes \EXT \ARR 
H^0(\omega_{\ZE} \otimes \OO_X (-L -K_X))= \left(\HOZK \right)^{\ast}
\end{equation}
which is equivalent to the sheaf morphism $\theta_{[e]}$ asserted in the lemma.
Turning to the second assertion of the lemma we recall the Koszul sequence (\ref{Koz})
and the extension class $\AL_e$ defined by it. The pairing (\ref{Ye}) with this class
gives a homomorphism
$$
Y_e (\cdot,\AL_e) : \HOZ \ARR \left(\HOZK \right)^{\ast}
$$
which is known to be an isomorphism given explicitly by the residue map as in
Chapter 5,\cite{[G-H]}.
\end{pf} 
\begin{cor}\label{div-e}
The locus where $\theta_{[e]}$ in (\ref{theta-e}) fails to be an isomorphism is a hypersurface
${\bf \Theta} ([e])$ of degree $d$ in $\PP(\EXT)$. Furthermore for $[e] \in \RUE$ and
$z \in Z_e$ the subset
\begin{equation}\label{th-z-e}
{\bf \Theta} ([e],z))= \{[\AL] \in \PP(\EXT) \mid \AL(z)= 0 \}
\end{equation}
is a hyperplane in $\PP(\EXT)$ and the support of ${\bf \Theta} ([e])$ is equal to
$\displaystyle{\bigcup_{z\in \ZE} {\bf \Theta} ([e],z))}$.
\end{cor}
\begin{pf}
Since the morphism $\theta_{[e]}$ in (\ref{theta-e}) is generically an isomorphism
its determinant 
$\det \theta_{[e]}$ is a nonzero section of 
$\OO_{\PP(\EXT)} (d)$ and the locus where $\theta_{[e]}$ fails to be an isomorphism is the
 hypersurface
$$
{\bf \Theta} ([e]) = (\det \theta_{[e]} =0)
$$
The second assertion of the corollary follows from the identification (\ref{SD4}) and the fact
that the dualizing sheaf $\omega_{\ZE}$ of $\ZE$ is an invertible sheaf on $Z_e$. This implies
that we can identify 
$\EXT$ with the subspace of sections of an invertible sheaf on $\ZE$. This and the fact that
there is a section in $\EXT$ which does not vanish anywhere on $\ZE$ imply that
${\bf \Theta} ([e],z))$ in (\ref{th-z-e}) is a hyperplane in 
$\PP(\EXT)$, for every $z\in Z_e$.
It is then obvious that the union of those hyperplanes is contained in ${\bf \Theta} ([e])$.
On the other hand it is clear that the complement of 
$\displaystyle{\bigcup_{z\in \ZE} {\bf \Theta} ([e],z))}$ in $\PP(\EXT)$ is contained in the
complement of ${\bf \Theta} ([e])$. Hence the equality
$$
\bigcup_{z\in \ZE} {\bf \Theta} ([e],z)) = support \,\,\,of\,\,\,{\bf \Theta} ([e])
$$
\end{pf}
As $[e]$ varies in $\UE^{\prime}$ the hypersurfaces ${\bf \Theta} ([e])$ will vary to form
a divisor in $\JUE$ (see \cite{[R1]},\S1.2, for technical details). This divisor will be
called the {\it theta-divisor} of $\JUE$ and denoted ${\bf \Theta}_{\UE}$.
\begin{defi}\label{ext-reg}
An extension class $\AL \in \EXT$ is called regular if the morphism
$\theta_{[e]}$ is an isomorphism at the point $[\AL] \in \PP(\EXT)$ or,
equivalently, if the homomorphism
$$
Y_ {[e]} (\cdot,\AL): \HOZ \ARR (\HOZK)^{\ast}
$$
induced by the pairing $Y_ {[e]}$ in (\ref{Ye}) is an isomorphism.
\end{defi}
Once we have the morphism $\theta_{[e]}$ in (\ref{theta-e}) we can define a
distinguished subsheaf of $\HOZ \otimes \OO_{\PP(\EXT)}$ (which will eventually become the first step of our filtration). This is done as follows.
Tensor the morphism $\theta_{[e]}$  with $\OO_{\PP(\EXT)} (1)$ to obtain the morphism
\begin{equation}\label{theta-e1}
\theta_{[e]} (1) : \HOZ \otimes \OO_{\PP(\EXT)} \ARR (\HOZK)^{\ast}\otimes \OO_{\PP(\EXT)} (1)
\end{equation}
Define the subsheaf $\HT_{[e]}$ of $\HOZ \otimes \OO_{\PP(\EXT)}$ to be the inverse image
 of the subsheaf $\EXT \otimes \OO_{\PP(\EXT)} (1)$ of 
$(\HOZK)^{\ast}\otimes \OO_{\PP(\EXT)} (1)$ under the morphism $\theta_{[e]} (1)$. 
The following properties of the sheaf $\HT_{[e]}$ follow immediately from its definition and Corollary \ref{div-e}.
\begin{pro}\label{pro-HTe}
\begin{enumerate}
\item[1)]
The sheaf $\HT_{[e]}$ contains $\OO_{\PP(\EXT)}$ as its subsheaf.
\item[2)]
The morphism $\theta_{[e]} (1)$ induces the following commutative diagram
\begin{equation}\label{theta-e-di}
\xymatrix{
0\ar[d]& 0\ar[d]\\
\OO_{\PP(\EXT)} \ar[d] \ar@{=}[r]& \OO_{\PP(\EXT)}\ar[d] \\
\HT_{[e]} \ar[d] \ar[r]^(.3){\theta_{[e]} (1)} & \EXT \otimes \OO_{\PP(\EXT)} (1) \ar[d] \\
\HT_{[e]} / \OO_{\PP(\EXT)} \ar[d] \ar[r]^{\theta^{\prime}_{[e]} (1)} & {\cal T}_{\PP(\EXT)} \ar[d]\\
0&0 }
\end{equation}
where ${\cal T}_{\PP(\EXT)}$ is the tangent bundle of $\PP(\EXT)$ and the column on the right
is the Euler sequence on $\PP(\EXT)$.
\item[3)]
The morphism $\theta_{[e]} (1)$ (resp. $\theta^{\prime}_{[e]} (1)$) in (\ref{theta-e-di})
is an isomorphism precisely on the complement of the divisor ${\bf \Theta}([e])$ in
Corollary \ref{div-e}. In particular, 
$\HT_{[e]}$ is a torsion free sheaf on $\PP(\EXT)$ 
 and its singularity set is ${\bf \Theta}([e])$.
Furthermore, for general $[e] \in \UE$, the rank of $\HT_{[e]}$ is equal to $\delta(\UE)$.
\end{enumerate}
\end{pro}
\begin{rem}
If $[\AL]$ is a point of $\PP(\EXT)$ corresponding to a regular extension class $\AL$, then
from the definition of the morphism $\theta_{[e]} $ in (\ref{theta-e}) it follows that the fibre $\HT_{[e]} ([\AL])$ of $\HT_{[e]}$ is the following subspace of functions on $\ZE$
\begin{equation}\label{HT-e-a}
\HT_{[e]} ([\AL]) = \{ \frac{\beta}{\alpha} \mid \beta \in \EXT \}
\end{equation}
where we think of elements of $\EXT$ as sections of the invertible sheaf
$\omega_{\ZE} \otimes \OO_X (-L-K_X)$
and the quotient of two sections of an an invertible sheaf is a (rational) function.
Furthermore the fact that $\AL$ is regular means that it is a section of
$\omega_{\ZE} \otimes \OO_X (-L-K_X)$ which does not vanish anywhere on $\ZE$.
Hence the quotient $\displaystyle{\frac{\beta}{\alpha}}$ is a (regular) function on $Z_e$.
\end{rem}
We are now in the position to define our filtration on
$\HOZ \otimes \OO_{\PP(\EXT)}$.
\\
\indent
1) Define
\begin{equation}\label{HT-1}
\HT_{-1} ([e]):=\HT_{[e]} 
\end{equation}
to be the first step of filtration. 
\\
\indent
2) To define its next steps we use the multiplication in
$\HOZ$. Namely, for $k\geq 2$ consider the morphism
\begin{equation}\label{mk}
m_k ([e]): S^k \HT_{[e]} \ARR \HOZ \otimes \OO_{\PP(\EXT)}
\end{equation}
induced by the multiplication in $\HOZ$. Define
\begin{equation}\label{HT-k}
\HT_{-k} ([e]):= im(m_k ([e]) )
\end{equation}
the image of $m_k ([e])$.
\\
\indent
The inclusion $\HT_{-k} ([e]) \subset \HT_{-k-1} ([e])$ follows from the fact that the constant
function $1$ is a global section of $\HT_{-1} ([e])$ (Proposition \ref{pro-HTe},1)).
Thus we obtain a decreasing filtration of $\HOZ \otimes \OO_{\PP(\EXT)}$
\begin{equation}\label{HTe-filt}
\HT_{-1} ([e])\subset \ldots \subset \HT_{-k} ([e]) \subset \HT_{-k-1} ([e]) \subset \ldots
\end{equation}
\begin{rem}\label{geomm}
Let $\HT_{-\bullet} \EA$ be the fibre of the filtration (\ref{HTe-filt}) at a point
$[\AL] \in \PP(\EXT)$. This filtration of $\HOZ$ has a simple geometric interpretation.
Namely, the subspace $\HT_{-1} \EA$ of functions on $Z_e$ defines the obvious morphism
\begin{equation}\label{kappaEA}
\tilde{\kappa} \EA : \ZE \ARR  \left(\HT_{-1} \EA \right)^{\ast}
\end{equation}
Then the subspace $\HT_{-k} \EA$ is the pullback of the homogeneous polynomial
functions of degree $k$ on $\left(\HT_{-1} \EA \right)^{\ast}$, i.e.
$$
\HT_{-k} \EA = \tilde{\kappa}^{\ast} \EA ( S^k \HT_{-1} \EA)
$$
In particular, the last piece of the filtration $\HT_{-\bullet} \EA$ 
is the subring
$\tilde{\kappa}^{\ast} (H^0(\OO_{Z^{\prime}_e}))$,
where $Z^{\prime}_e$ is the image of $\tilde{\kappa} \EA$.
\end{rem}
We have described our filtration for a fixed $[e]$. As $[e]$ varies in $\UEP$ the
sheaves in (\ref{HTe-filt}) fit together to give the decreasing filtration
\begin{equation}\label{HT-filt}
\HT_{-1} \subset \ldots \subset \HT_{-k}  \subset \HT_{-k-1}  \subset \ldots
\end{equation}
of the locally free sheaf $\PIP$. In particular, all sheaves in the filtration are torsion free since they are nonzero subsheaves of a locally free sheaf. Hence the ranks of the
sheaves in  (\ref{HT-filt}) are well-defined.
\begin{defi}\label{HilbUE}
Define the function $P_{\UE} (k) = rk \HT_{-k}$, for every integer $k \geq 1$.
This function  is called the Hilbert function of the universal cluster of $\EE$ over $\UE$.
The smallest value $w(\UE)$ such that
$$
P_{\UE} (k) =P_{\UE} (w(\UE)), \forall k\geq w(\UE)
$$
is called the weight of of the universal cluster of $\EE$ over $\UE$.
\end{defi}
\begin{rem}
\begin{enumerate}
\item[1)]
By definition and Proposition \ref{pro-HTe} 
$$
P_{\UE} (1) = \delta(\UE)
$$
\item[2)]
$P_{\UE} (w(\UE)) = deg Z^{\prime}_e $, where $Z^{\prime}_e $ is the image of the morphism
$\tilde{\kappa} \EA$ in  (\ref{kappaEA}), for general $\EA \in \JUE$.
\end{enumerate}
\end{rem}
 Let $S$ be the union of the singularity sets of the sheaves
$\HT_{-k},k=1,\ldots,w(\UE)$, and let $\JUE^{\prime}$ be the complement of $S$
in $\JUE$. This is the largest Zariski open subset of $\JUE$ where the filtration
(\ref{HT-filt}) is locally free, i.e. every sheaf in (\ref{HT-filt}) is locally free over
$\JUE^{\prime}$.
\begin{rem}
Let $\UE^{\prime \prime} = \pi_{\UE} (\JUE^{\prime})$. This is a nonempty Zariski open subset of $\UE$. Let ${\bf \Theta}^{\prime \prime}_{\UE}$ be the part of the theta-divisor
${\bf \Theta}_{\UE}$ lying over $\UE^{\prime \prime}$, i.e.
$$
{\bf \Theta}^{\prime \prime}_{\UE} ={\bf \Theta}_{\UE} \cap \pi^{-1}_{\UE} (\UE^{\prime \prime})
$$
Then we claim that
\begin{equation}\label{JUp}
\JUE^{\prime} = \pi^{-1}_{\UE} (\UE^{\prime \prime}) \setminus 
{\bf \Theta}^{\prime \prime}_{\UE}
\end{equation}
This follows from the fact that for every $[e]\in \UE$
the sheaves of the filtration (\ref{HTe-filt}) have constant ranks on the complement of the divisor ${\bf \Theta}([e])$ in Corollary \ref{div-e}.
Indeed, let $[\AL]$ and $[\AL^{\prime}]$ be two distinct points of
$\PP(\EXT) \setminus {\bf \Theta}([e])$. Then from (\ref{HT-e-a}) it follows that the
multiplication by the function 
$\displaystyle{f=\frac{\AL}{\AL^{\prime}}}$ induces an isomorphism
$$
f: \HT_{[e]} ([\AL]) \ARR \HT_{[e]} ([\AL^{\prime}])
$$
This and the definition of $\HT_{-k} ([e])$ implies that the multiplication
$$
f^k : \HT_{-k}\EA \ARR \HT_{-k} ([e],[\AL^{\prime}])
$$
is an isomorphism between the fibres of 
$\HT_{-k}$
at $\EA$ and $([e],[\AL^{\prime}])$, for every $k\geq 1$.
\end{rem}

\section{ Theorem \ref{Th1} and Theorem \ref{Th2}:
mathematics}
In this section we give proofs (or their sketch with reference to \cite{[R1]} for details) of the results described in
{\bf Part II} of the introduction. 
\subsection{Mathematics of Theorem \ref{Th1}}
 
The direct sum decomposition (\ref{dsds}) in Theorem \ref{Th1} is a result of splitting of the filtration (\ref{HT-filt}). To explain this we need to observe that over the locus
$\RUE$ parametrizing  regular sections of $\EE$ the sheaf
$p_{\RUE \ast} \OO_{\ZDO}$ acquires an extra-structure of being a sheaf of Frobenius
algebras. Indeed, for $[e] \in \RUE$, the ring $\HOZ$ comes equipped with the trace map
$$
Tr_{[e]} : \HOZ \ARR {\CC}
$$
defined by the formula
\begin{equation}\label{tr-e}
  Tr_{[e]} (f) = \sum_{z\in \ZE} f(z),\,\,\,\forall f \in \HOZ
\end{equation}
The trace map gives rise to a non-degenerate symmetric bilinear pairing
\begin{equation}\label{q-e}
q_{[e]} : \HOZ \times \HOZ \ARR {\CC}
\end{equation}
defined as follows
\begin{equation}\label{q-e1}
q_{[e]} (f,g) = Tr_{[e]} (fg), \,\,\,\forall f,g \in \HOZ
\end{equation}
If no confusion is likely we will not distinguish the above pairing and the quadratic form associated to it and refer to both as quadratic form.
\\
\\
\indent
 The essential observation is as follows.
\begin{lem}\label{q-H-k}
There is a nonempty Zariski open subset $\JBBU$ of $\JUE^{\prime}$ lying over
$\RUE^{\prime\prime} =\RUE \cap \UE^{\prime\prime}$ such that for every
$\EA \in \JBBU$ the restriction of the quadratic form
$q_{[e]}$ to the subspaces $\HT_{-k} \EA,\,\,\forall k\geq1$, is again non-degenerate.
\end{lem}
For a proof of this result we refer to \cite{[R1]},\S2,Claim 2.2.
\\
\\
\indent
The quadratic forms $q_{[e]}$ fit together to define the quadratic form
$\QB$ on the sheaf $\PIPO$. Using this form we define
\begin{equation}\label{HBp}
\HB^p = \left\{
\begin{array}{cl}
\left(\HT_{-p}\right)^{\bot} \cap \HT_{-p-1},&if\,\,\,p=0,\ldots,w(\UE)-1\\
\left(\HT_{-w(\UE)}\right)^{\bot},&if\,\,\,p=w(\UE)
\end{array}
\right.
\end{equation}
where $(\HT_{-s})^{\bot}$ denotes the orthogonal complement of
$\HT_{-s}$ in $\PIPO$ with respect to the quadratic form $\QB$ and where we agree that
$\HT_{0} =0$.
\\
\indent
From Lemma \ref{q-H-k} it follows that over the open set $\JBBU$ we have an orthogonal
direct sum decomposition
\begin{equation}\label{orthd}
\PIPO \otimes \OO_{\JBBU} = \bigoplus^{w(\UE)}_{p=0} \HB^p_{\JBBU}
\end{equation}
where $ \HB^p_{\JBBU} = \HB^p \otimes \OO_{\JBBU}$. This gives us the first statement 
of Theorem \ref{Th1}.
\\
\indent
From now on, unless said otherwise,
we will be working over $\JBBU$ and we omit the subscripts from the notation above,
if no confusion is likely.
\\
\\
\indent
The summand $\HB^0$ plays a special role in our constructions. First of all, by definition (\ref{HBp}),
we have the following
\begin{equation}\label{H0}
\HB^0 = \HT_{-1} \otimes \OO_{\JBBU}
\end{equation}
This identification implies the relative version of Proposition \ref{pro-HTe} which
is nothing but properties $(a)$ and $(b)$ of Theorem \ref{Th1},2).
\begin{defi}\label{ext-p}
Let $[e] \in \RUE$. A regular extension class $\AL \in \EXT$ (see Definition \ref{ext-reg})
is called polarizing or a polarization of $\ZE$ if the quadratic form
$q_{[e]}$ in (\ref{q-e}) is non-degenerate on $\HT_{-k} \EA$, for every $k\geq 1$.
\end{defi}
Our terminology is of course borrowed from the usual Hodge decomposition of the cohomology
of compact K\"{a}hler manifolds induced by a K\"{a}hler metric or a polarization on a manifold.
In fact the decomposition (\ref{orthd}) could be regarded as an analogue of Hodge decomposition. With this analogy in mind a polarizing extension class $\AL \in \EXT$
could be envisaged as a kind of non-classical ``K\"{a}hler" structure on $Z_e$. Thus the
vertical (fibre) directions of the variety $\JBBU$ could be viewed as the moduli space
of non-classical ``K\"{a}hler" structures, while in the horizontal directions one has the classical variation of configurations $\ZE$ as $[e]$ varies in
$\RUE^{\prime\prime}$.

\subsection{Mathematics of Theorem \ref{Th2}}
The direct sum decomposition (\ref{orthd}) can be viewed as a variation of the decomposition
\begin{equation}\label{orthdEA}
\HOZ = \bigoplus^{w(\UE)}_{p=0} \HB^p \EA
\end{equation}
as a point $\EA$ moves in $\JBBU$ and where 
$\HB^p \EA$ is the fibre of $\HB^p_{\JBBU}$ at $\EA$. The key observation is that the variation
with respect to the parameter $[\AL]$ is entirely captured by the multiplication in the ring
$\HOZ$ by functions lying in the summand $\HB^0 \EA$.
This follows from the identification (\ref{H0}) and the isomorphism 
$\theta^{\prime}$ in (\ref{theta-d}).
So we will concentrate on the properties of the multiplication by functions in 
$\HB^0 \EA$. In particular, we want to understand how this multiplication acts on the
direct sum decomposition (\ref{orthdEA}).
\begin{lem}\label{act-d-t}
Let
$D(t): \HOZ \ARR \HOZ$ be the homomorphism of multiplication by
$t\in \HB^0 \EA$.
Then the following holds.
\begin{enumerate}
\item[1)]
The subspace
\begin{equation}\label{H-w-EA}
\HT_{-w(\UE)} \EA = \bigoplus^{w(\UE)-1}_{p=0} \HB^p \EA
\end{equation}
is preserved by $D(t)$.
\item[2)]
Let $D_p (t)$ be the restriction of $D(t)$ to the summand
$\HB^p \EA$. Then
\begin{equation}\label{Dp-t-act}
D_p (t): \HB^p \EA \ARR \HB^{p-1} \EA \oplus \HB^p \EA \oplus \HB^{p+1} \EA
\end{equation}
for every $p=1,\ldots,w(\UE)-2$. Furthermore, for $p=0,w(\UE)-1$, the homomorphism $D_p (t)$ has the following form
\begin{equation}\label{D0-t-act}
D_0 (t): \HB^0 \EA \ARR  \HB^0 \EA \oplus \HB^{1} \EA
\end{equation}
\begin{equation}\label{Dw-t-act}
D_{w-1} (t): \HB^{w-1} \EA \ARR  \HB^{w-1} \EA \oplus \HB^{w-2} \EA
\end{equation}
\end{enumerate}
\end{lem}
\begin{pf}
The first statement follows from the definition of $w(\UE)$. To prove the second we take
$h \in \HB^p \EA = \left(\HT_{-p} \EA\right)^{\bot} \cap \HT_{-p-1}\EA$. Multiplying it by
$t\in \HB^0 \EA$ we obtain
\begin{equation}\label{th}
D_p (t)(h)=th \in \HT_{-p-2}\EA =\HT_{-p}\EA \oplus \HB^p \EA \oplus \HB^{p+1} \EA
\end{equation}
Furthermore since $h$ is also in $\left(\HT_{-p} \EA\right)^{\bot}$ we have
$$
q_{[e]}(th,h^{\prime}) =q_{[e]}(h,th^{\prime}) =0,\,\,\,\forall h^{\prime} \in \HT_{-p+1}\EA
$$
This implies that $th \in \left(\HT_{-p+1}\EA\right)^{\bot}$. Combining with
(\ref{th}) yields
\begin{eqnarray*}
D_p (t)(h)=&th \in \left(\HT_{-p+1}\EA\right)^{\bot} \cap \left(
\HT_{-p}\EA \oplus \HB^p \EA \oplus \HB^{p+1} \EA \right)\\
 &=
\HB^{p-1} \EA \oplus \HB^p \EA \oplus \HB^{p+1} \EA
\end{eqnarray*}
\end{pf}
Denote by 
$D^{-}_p (t),D^{0}_p (t),D^{+}_p (t)$ the projections of 
$D_p (t)$ on the first, second and the third summand of
$\HB^{p-1} \EA \oplus \HB^p \EA \oplus \HB^{p+1} \EA$ respectively.
Set 
\begin{equation}\label{Dt+-0}
D^{\pm}(t) = \sum^{w-1}_{p=0} D^{\pm}_p(t)\,\,\,and\,\,\,
D^{0}(t) = \sum^{w-1}_{p=0} D^{0}_p(t)
\end{equation}
Then the multiplication homomorphism $D(t)$, restricted to the subspace
 $\HT_{-w} \EA$, admits the following decomposition
\begin{equation}\label{dec-Dt}
D(t) =D^{-}(t) + D^{0}(t) + D^{+}(t)
\end{equation}
The components $D^{\pm}(t)$ act on the decomposition in (\ref{H-w-EA}) as homomorphisms
of degree $\pm1$, i.e. they shift the grading in (\ref{H-w-EA}) by $\pm1$,
while $D^{0}(t)$ is a homomorphism of degree $0$, i.e. it preserves the decomposition
in (\ref{H-w-EA}).
Thus with every point $\EA \in \JBBU$ we can associate the subspace
\begin{equation}\label{sp-end}
\{ D^{0}(t),D^{\pm}(t) \mid t\in \HB^0 \EA \}
\end{equation}
of endomorphisms of $\HT_{-w} \EA$. These endomorphisms could be viewed as rules for the infinitesimal variation of the direct sum
 decomposition in (\ref{H-w-EA}) along the directions of the fibre of $\JBBU$ at $\EA$.
\\
\indent
The family of the endomorphisms in (\ref{sp-end}) can be rewritten as homomorphisms
\begin{equation}\label{D-EA}
D^{0} \EA,\,\,D^{\pm} \EA : \HT_{-w} \EA \ARR 
\left(\HB^0 \EA \right)^{\ast} \otimes \HT_{-w}\EA
\end{equation}
Finally, letting $\EA$ vary in $\JBBU$ the above homomorphisms fit together to yield the
morphisms of sheaves
\begin{equation}\label{D0pm}
D^{0},D^{\pm} : \HT_{-w}  \ARR \left(\HB^0  \right)^{\ast} \otimes \HT_{-w}
\end{equation}
while the multiplication morphism
$$
D: \HT_{-w}  \ARR \left(\HB^0  \right)^{\ast} \otimes \HT_{-w}
$$
admits the following decomposition
\begin{equation}\label{mor-dec}
D= D^{-} +D^{0} +D^{+}
\end{equation}
The commutativity of multiplication now translates into the fact that
$D$ is a Higgs morphism of $\HT_{-w}$ with values in $\left(\HB^0  \right)^{\ast}$
(see \cite{[S]} or Definition 3.2,\cite{[R1]}), i.e.
\begin{equation}\label{com}
D^2 = D \wedge D =0
\end{equation}
This together with the decomposition in (\ref{mor-dec}) yields
$$
0 = (D^{-} +D^{0} +D^{+} ) \wedge (D^{-} +D^{0} +D^{+} )
$$
Expanding and grouping together the terms of the same degree yield the relations
(\ref{rel}) of Theorem \ref{Th2}.

\subsection{Trivalent graph associated to $\JBBU$}
The data of the direct sum decomposition (\ref{HT-wS}) and the morphisms 
$D^{0},D^{\pm}$ of Theorem \ref{Th2} can be organized in the following trivalent graph
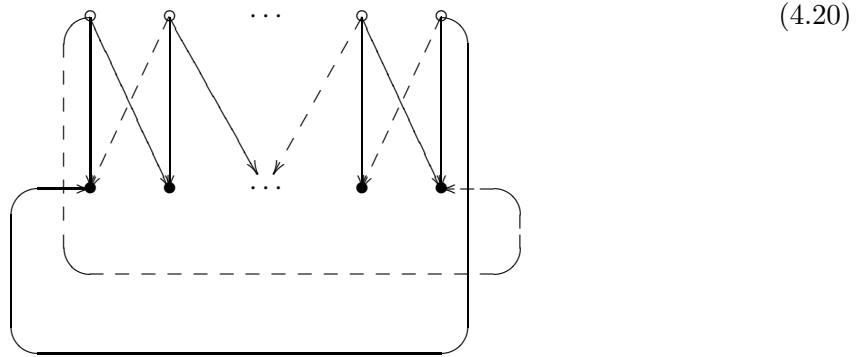
\begin{equation}\label{G}
\xymatrix{
&*=0{\circ} \ar[2,0] \ar[2,1] \ar@{--{>}}`l[d]                               
                               `[3,1]                               
                               `[2,5]                               
                                `[2,4]
                                [2,4]
                              &
                               *=0{\circ}\ar[2,0] \ar[rdd] \ar@{--{>}}[ddl]&
\cdots&*=0{\circ} \ar[2,0] \ar@{--{>}}[ldd] \ar[2,1]&
*=0{\circ} \ar[2,0] \ar@{--{>}}[2,-1] \ar`r[d] 
                              `[4,-2]
                              `[2,-5]                              
                              `[2,-4]
                              [2,-4]&\\
& & & & & & \\
&*=0{\bullet}&*=0{ \bullet}&\cdots&*=0{\bullet}&*=0{\bullet}&\\
& & & & & &\\
& & & & & &}
\end{equation}

There are $w=w(\UE)$ pairs of vertices in (\ref{G}) which are ordered by the index set
$I= \{0,\ldots,w-1\}$ from left to right with the $i$-th pair being the ends of the 
$i$-th vertical edge in (\ref{G}). This vertices are labeled by $i_u$ and $i_d$, for
the upper  and lower one, respectively. Besides vertical edges  relating
$i_u$ with $i_d$, for $i\in I$, the other edges are drawn as follows.
For every $i\in I$ the vertex $i_u$ is connected to $(i-1)_d$ and $(i+1)_d$
with understanding that $(-1)_d = (w-1)_d$ and $w_d = 0_d$. This graph 
will be denoted by $G(\JBBU)$
or simply by $G$ if no ambiguity is likely.
\\
\indent
We agree that all edges of $G$ are oriented from ``up" to ``down".
The graph $G$ with this orientation will be denoted by $\vec{G}$. In this 
orientation the edges incident to an upper (resp. lower) vertex of $G$ are
always out-going (resp. in-going). At every upper vertex of $G$ we fix a 
counterclockwise
order of incident edges and color them, starting with the vertical one,
by ``$0$", ``$+$", and ``$-$", respectively. This will be called the 
{\it natural coloring of } $G$. 
\\
\indent
It should be clear that the graph $G$ with its coloring and orientation is just a pictorial
representation of the decomposition
(\ref{HT-wS}) and the morphisms $D^0,D^{\pm}$ in Theorem \ref{Th2}.
Namely, if we place the summand $\HB^i$ of the direct sum decomposition at the vertices
$i_u$ and $i_d$, for every $i \in I$, and label the (oriented) edges at $i_u$
 by morphisms $D^0,D^{\pm}$ so that the natural
colors of edges are respected, i.e. we place $D^0$ on the edge colored by
``0" and $D^{\pm}$ on the edges colored by ``$\pm$" respectively, then we will
recover the rules of action of the morphisms 
$D^0,D^{\pm}$ on the summand $\HB^i$ placed at the vertex $i_u$.
We will see shortly how one can use this pictorial representation to obtain ``quantum"
properties of the universal configuration.

\section{Theorem \ref{Th3} and Theorem \ref{Th4}}
The results of Theorem \ref{Th3} are obtained by reinterpreting the morphisms
$D,D^0,D^{\pm}$ in Theorem \ref{Th2} in geometric terms, while 
Theorem \ref{Th4} is an application of path-operators (discussed in detail in \S6) in the 
geometric context of Theorem \ref{Th3}. In this section it always will be assumed that the
weight $w=w(\UE) \geq 2$.
\subsection{Theorem \ref{Th3}}
The main idea is to take the morphism $D$ in (\ref{D}) and to perturb it according to the
direct sum decomposition (\ref{HT-wS}) in a such a way that the commutativity condition in
(\ref{com})
would be preserved. More precisely, set 
$D^0_p$ (resp. $D^{\pm}_p$) to be restrictions of the morphism $D^0$ (resp. $D^{\pm}$) to
the $p$-th summand of the decomposition in (\ref{HT-wS}). Then $D$ admits the following
decomposition
\begin{equation}
D = D^0 + D^{+} + D^{-} = \sum^{w-1}_{p=0} D^0_p + \sum^{w-2}_{p=0} D^{+}_p + 
\sum^{w-2}_{p=0} D^{-}_{p+1}
\end{equation}
We take now a sufficiently general deformation of $D$ of the form
\begin{equation}\label{txy}
D({\bf t},{\bf x},{\bf y}) = \sum^{w-1}_{p=0} t_p D^0_p  + \,\,
 \sum^{w-2}_{p=0} x_p D^{+}_p\,\,
+ \sum^{w-2}_{p=0} y_{p} D^{-}_{p+1}
\end{equation}
where 
${\bf t}=( t_p) \in {\bf C^w},\,\,{\bf x}=(x_p),\,\,{\bf y}=(y_{p}) \in {\bf C^{w-1}} $,
and we seek conditions on the parameters ${\bf t},{\bf x},{\bf y}$ such that
 \begin{equation}\label{eq-Higgs}
D^2 ({\bf t},{\bf x},{\bf y}) = D({\bf t},{\bf x},{\bf y}) \wedge D({\bf t},{\bf x},{\bf y})=0
\end{equation}
would hold, i.e. $D({\bf t},{\bf x},{\bf y})$ is a Higgs morphism of $\HT_{-w}$ with values in
$(\HB^0)^{\ast}$. 
\\
\indent
Using the relations (\ref{rel}) it is straightforward 
 to deduce the following (see \cite{[R1]},\S4, for details).
\begin{lem}\label{H-hat} 
The set
\begin{equation}\label{sol}
\hat{H}_{\JBBU} = \{ (z,{\bf x},{\bf y})\in {\bf C\times C^{w-1}\times C^{w-1}} \mid
{\bf x}=(x_p),\,\,{\bf y}=(y_p),\,\, x_p y_p =z^2,\,\,p=0,\ldots,w-2 \}
\end{equation}
is a set of solutions for (\ref{eq-Higgs}), i.e. the morphism
$$
D(z,{\bf x},{\bf y}) = zD^0 + \sum^{w-2}_{p=0} x_p D^{+}_p + 
\sum^{w-2}_{p=0} y_p D^{-}_{p+1}
$$
satisfies (\ref{eq-Higgs}), for all $(z,{\bf x},{\bf y}) \in \hat{H}_{\JBBU}$.
\end{lem}
The set $\hat{H}_{\JBBU}$ is an affine variety. It admits an obvious $\bf C^{\ast}$-action, since scaling a solution of (\ref{eq-Higgs}) by a constant is again a solution 
of (\ref{eq-Higgs}).
\begin{defi}\label{Alb}
The set 
$$
H_{\JBBU} = \hat{H}_{\JBBU} / {\bf C^{\ast}}
$$
of homothety equivalent Higgs morphisms parametrized by $\hat{H}_{\JBBU}$ is called the nonabelian Albanese of $\JBBU$.
\end{defi}
As usual, if no ambiguity is likely, we will omit the subscript in $H_{\JBBU}$.
\\
\indent 
We will now give a projective description of the Albanese $H$. 
\begin{pro}\label{H1}
Let $V$ be the vector space freely generated by symbols
$v_0,v^{\pm}_p,(p=0,\ldots,w-2)$
and let $V^{\ast}$ be the dual of $V$ with the dual basis
$T,X_p,Y_p,\,(p=0,\ldots,w-2)$. Then
$H_{\JBBU}$ is a subvariety of $\PP(V)$ defined by the quadratic equations
$$
 X_p Y_p =T^2 \,(p=0,\ldots,w-2)
$$
In particular, $H_{\JBBU}$ is a projective variety of dimension $w-1$ and it has degree
$2^{w-1}$ in $\PP(V)$. Furthermore its dualizing
sheaf $\omega_H = \OO_H (-1)$, i.e. $H_{\JBBU}$ is a Fano variety of dimension
$w-1$. 
\end{pro}
This projective description of $H_{\JBBU}$ follows immediately from the description of 
$\hat{H}_{\JBBU}$ in Lemma \ref{H-hat}.
\begin{cor}\label{h-sec}
Let $w=w(\UE) \geq 3$. Then hyperplane sections of $H$ are Calabi-Yau varieties of dimension
$w-2$.
\end{cor}
\begin{pf}
The statement about the dualizing sheaf $\omega_H$ of $H$ in
Proposition \ref{H1} and 
the adjunction formula (see e.g. \cite{[G-H]}) yield the assertion.
\end{pf}
For other properties of $H$, in particular, its toric description we refer to \S4,\cite{[R1]}.

\subsection{Theorem \ref{Th4}}\label{sCY}
This theorem establishes a geometric correspondence between the nonabelian Jacobian
$\JBBU$ and its nonabelian Albanese $H_{\JBBU}$. This correspondence is given by the morphism
\begin{equation}\label{CY-map}
CY_{\RUE}: \JBBU \ARR \PP(H^0(\OO_{H_{\JBBU}} (d))
\end{equation}
which is called the Calabi-Yau cycle map for an obvious reason that the divisors assigned to
points of $\JBB$ are cycles of hyperplane sections of $H$ which are, by Corollary \ref{h-sec},
Calabi-Yau varieties. Before giving a construction of this map let us consider the situation
from a more conceptual point of view.
\\
\indent
We have constructed two varieties naturally associated to the universal configuration
$\ZDO$. The first one is $\JBBU$ which as we have seen can be thought of as a moduli space of ``classical" an ``hidden" parameters of the universal configuration. The second one is
the nonabelian Albanese $H_{\JBBU}$ which is a moduli of homothety equivalent Higgs structures on the sheaf $\HT_{-w} \otimes \OO_{\JBBU}$. It should be observed that 
$H_{\JBBU}$ depends only on the weight $w=w(\UE)$ and the relations (\ref{rel}) between
morphisms $D^0,D^{\pm}$. In particular, it completely ``forgets" about our
variety $X$ and the universal configuration $\ZDO$. So a correspondence between
$\JBBU$ and its nonabelian Albanese $H$ could be viewed as some kind of duality which
partially restores ties of $H$ with the variety $X$. 
\\
\\
\indent
We now turn to the construction of the Calabi-Yau cycle map $CY_{\RUE}$. It is the same
as in the proof of Proposition 5.1,\cite{[R1]}, but we reproduce it here as an 
illustration/application of the path-operator technique (discussed in details in \S6).
\\
\indent
Let $\EA$ be a point of $\JBB$. To define 
$CY \EA$ we will produce a divisor of the form
$$
\sum_{z\in Z_e} H_{z,[\AL]}
$$
where $H_{z,[\AL]}$ is a hyperplane section of $H_{\JBB}$ intrinsically associated to
$(z,[\AL]) \in (Z_e,[\AL])$.
  To do this we recall
that $H^0(\OO_H (1))$ comes with a particular basis
$T,X_p,Y_p \,(p=0,\ldots,w-2)$ (Proposition \ref{H1}). So our strategy is as follows.
For every point $z\in Z$ we will produce the 
constants $t(z,[\alpha]),x_p(z,[\alpha]), y_p(z,[\alpha]) (p=0,\ldots,w-2)$
depending holomorphically on $(z,[\alpha])$. Then we use these constants 
to define 
the section
\begin{equation}\label{s(z)}
s(z,[\alpha]) = t(z,[\alpha])T + \sum^{w-2}_{p=0}x_p(z,[\alpha])X_p
+ \sum^{w-2}_{p=0}y_p(z,[\alpha])Y_p \in H^0(\OO_H (1))
\end{equation}
Then we define 
$CY \EA$ to be the product
$\displaystyle{\prod_{z\in Z_e}s(z,[\alpha])}$. Thus our argument comes down to defining the constants
$t(z,[\alpha]),x_p(z,[\alpha]), y_p(z,[\alpha]) (p=0,\ldots,w-2)$.
This is done by using the orthogonal direct sum decomposition (\ref{HT-wS}) at
the point $\EA$
\begin{equation}\label{Orth-HodgeZ}
 \HT_{-w} \EA = \bigoplus^{w-1}_{p=0} {\HB^p} \EA
\end{equation}
and the morphisms $D^0,D^{\pm}$ in Theorem \ref{Th2}.
 \\
 \indent
 We proceed by picking a distinguished vector in ${\HB^0}\EA$ for 
$(z,[\AL])$. Namely, let $\delta_z$ be the delta-function on $Z$ supported at $z$,
 i.e. $\delta_z (z) =1$ and vanishes at other points of $Z$.
 Denote by $\delta^0_z$ the component of $\delta_z$ in
 ${\HB^0}\EA$. This is our distinguished vector alluded to above.
\\
\indent
Next we take the operator
 $D^{+}(\delta^0_z)$ and apply it repeatedly to $\delta^0_z$ to obtain (the right) moving string of functions
 $\delta^{(p)}_z = (D^{+}(\delta^0_z))^p(\delta^0_z)$, for $p=0,\ldots,w-1$, i.e. we move
$\delta^0_z$ along the shortest path of the graph $G(\JBBU)$ in (\ref{G}) connecting the
first vertical edge of the graph with the last one and passing through every vertical edge in between.
\\
\indent
 Once we arrive to $\delta^{(w-1)}_z \in {\HB^{w-1}} \EA$
 we apply the operator
 $D^{-}(\delta^0_z)$ to create (the left) moving string of functions
 $\delta^{(w-1),(w-1-m)}_z = (D^{-}(\delta^0_z))^m (\delta^{(w-1)}_z)$,
 for $m=0,\ldots,w-1$ (here we move $\delta^{(w-1)}_z$ back to the first vertical edge of the graph by retracing the shortest path used above).
 The desired constants are obtained essentially by evaluating
 all these functions at $z$. More precisely, define
 $$
 x_p(z,[\alpha]) =exp(\delta^{(p)}_z(z)),\,
 y_p(z,[\alpha]) =exp(\delta^{(w-1),(p+1)}_z (z)),\,
 t(z,[\alpha]) =exp(\delta^{(w-1),(0)}_z (z)),
 $$
 where $p=0,\ldots,w-2$.
 \\
 \indent
 The above construction depends only on the orthogonal decomposition
 (\ref{Orth-HodgeZ}). Since the latter varies holomorphically with
 $\EA$ we obtain that $CY \EA$ depends holomorphically
 on $\EA$ as well. This completes the proof of Theorem \ref{Th4}.

\section{String theoretic analogies}
Until now our construction of the variety $\JUE$ was motivated by purely algebro-geometric considerations. However, the higher $Ext$ groups,
 the essential ingredient of the
construction, have become familiar objects to physicists in connection with recent
advances in understanding of the role and nature of $D$-branes in string theory
(see \cite{[D]} and the references therein). We will now attempt to review the variety
$\JUE$ and its properties contained in Theorem \ref{Th1}, Theorem \ref{Th2} and Theorem \ref{Th4}
from the string theoretic perspective.
\\
\\
\indent
The point of departure of our construction is the Zariski open subset $\UE$, parametrizing zero-dimensional subschemes of $X$,
which are zero-loci of ``good" sections of $\EE$. Furthermore, the smaller Zariski open subset $\RUE$ parametrizes the configurations of points on $X$.
This could be viewed as the ``classical" space of configurations of points (particles) on $X$ associated to our vector bundle $\EE$. 
 The main idea of our construction is to consider
$0$-dimensional subschemes $Z_e$, where $[e]$ varies in $\UE$ together with all nonzero extension classes
$\AL \in \EXT$. This is a kind of thickening of $\UE$. We suggest to to view the resulting space $\JUE$ as a ``quantum" space  associated to $\EE$ and
 the Zariski open subset $\JBBU$ 
 defined in Theorem \ref{Th2} as the space of  ``quantum" configurations associated to $\EE$. Let us 
explain our thinking.
\\
\indent
Going from $\RUE$ to $\JUEO$ consists of considering a configuration of points $Z_e \,([e] \in \RUE)$ 
together with all possible (projectivized) non-zero extension classes 
$[\alpha] \in \PP(\EXT)$. One should think of $[\alpha]$ as an equivalence class of complexes of sheaves on 
$X$ (such as the Koszul complex in (\ref{Koz})) ``wrapped" on $Z_e$. So the pair
$(\ZE,[\alpha])$
could be envisaged as a $D$-brane $\ZE$ (of dimension zero) with the ``string"  $[\alpha]$= complex of sheaves on $X$
corresponding to   $[\alpha]$, attached to it. So the fibre of $\JUEO$ over $[e]$ could be viewed as a space of ``stringy"
structures one can put on $Z_e$.
\\
\indent
These  ``stringy" structures acquire a precise meaning for polarizing extension classes in the sense of Definition  \ref{ext-p}. For a polarizing
extension class $\alpha$ in $\EXT$ we have an explicit identification of the group $\EXT$ (``strings wrapped" on $\ZE$)
 with the subspace $\HT_{[e]} ([\alpha]) = \HB^0 \EA$
provided by the morphism $\theta_{[e]} (1)$ (see Proposition \ref{pro-HTe}, 2)) evaluated at $[\alpha]$:
\begin{equation}\label{tea}
\theta_{[e]} ([\alpha]): \HT_{[e]} ([\alpha]) = \HB^0 \EA \longrightarrow \EXT
\end{equation}
The inverse $\theta^{-1}_{[e]} ([\alpha])$ of this isomorphism can be viewed as a ``regular linearization" of ``strings"
 having ends on the $D$-brane
$\ZE$: it identifies ``strings" in $\EXT$ as a subspace of functions on $\ZE$ (this is the meaning of ``regular" here).
Furthermore, this subspace is the summand of degree $0$ in the direct sum decomposition
 \begin{equation}\label{dec-HT-wEA}
\HT_{-w} \EA = \bigoplus^{w-1}_{p=0} \HB^p \EA
\end{equation}
which we have seen in (\ref{H-w-EA}).
From (\ref{sp-end}) we also know that the functions in $\HB^0 \EA$ define the space of endomorphisms 
\begin{equation}\label{sp-end-s}
\{ D^0(t),D^{\pm} (t) \mid t \in \HB^0 \EA \}
\end{equation}
Putting (\ref{tea}) and (\ref{sp-end-s}) together we obtain linear (regular) representation of $\EXT$:
\begin{equation}
\EXT \longrightarrow \bigoplus^3 End (\HT_{-w} \EA)
\end{equation}
which assigns to a ``string" $\beta \in \EXT$ a triple of endomorphisms
\begin{equation}\label{autr}
D^0 (t_{\beta}), D^{\pm} (t_{\beta}) : \HT_{-w} \EA \longrightarrow  \HT_{-w}
\end{equation}
where $t_{\beta} = \theta_{[e]}^{-1} ([\alpha]) (\beta)$ is the function corresponding to $\beta$ under the isomorphism
in (\ref{tea}).
\\
\indent
The above considerations constitute quantum aspects of the space $\JBBU$.
Indeed, by definition $\JBBU$ lies over the space $\RUE^{\prime\prime}$ (see Lemma \ref{q-H-k}). This latter space can be viewed as 
the classical level parametrizing 
configurations $\ZE$, with $[e] \in \RUE^{\prime\prime}$.  Their spaces of functions $\HOZ$ can be viewed as classical observables on $\ZE$.
By going from $[e] \in \RUE^{\prime\prime}$ to $\JBBU$ reveals the parameter $\AL$ varying in the 
  fibre $\JBBU ([e])$ of $\JBBU$ over $[e]$.
Once the pair $\EA \in \JBBU$ is fixed the classical observables
 assigned to  ``strings" in $\EXT$, i.e. functions  in $\HB^0 \EA$, are transformed into
quantum observables = linear operators on $ \HT_{-w}$ in (\ref{autr}).  
Furthermore, the trivalent graph $G$ in (\ref{G}) can now be used
as an explicit mechanism to create many more such operators as well as to``propagate" our (functional representations of)
``strings" in the space $ \HT_{-w} \EA$. 
Namely, given a path $\gamma$ of the graph $G$ we define the path-operator
\begin{equation}\label{path-op}
D_{\gamma} \EA : \left(\HB^0 \EA \right)^{\otimes l(\gamma)} \ARR End(\HT_{-w} \EA)
\end{equation}
where $l(\gamma)$ is the length of $\gamma$, i.e. the number of edges composing
$\gamma$. This is done as follows.
\\
\indent
By definition $\gamma$ is a sequence $\{\gamma_1,\ldots,\gamma_{l(\gamma)} \}$
of oriented edges of $G$ such that the end of $\gamma_i$ is the origin of  $\gamma_{i+1}$,
for every $i=1,\ldots,l(\gamma)-1$. Define
\begin{equation}\label{path-op1}
D_{\gamma} \EA (t_1 \otimes \cdots \otimes t_{l(\gamma)})=
D^{c(\gamma_{l(\gamma)})}(t_{l(\gamma)}) \circ \cdots \circ D^{c(\gamma_1)}(t_1)
\end{equation}
for all $t_1,\ldots,t_{l(\gamma)} \in \HB^0 \EA $, where the right-hand side is the composition
of the endomorphisms in (\ref{sp-end-s}) with $c(\gamma_i) \in \{-,0,+ \}$ determined by the following rule:
\\
if the orientation of the edge $\gamma_i$ in $\gamma$ coincides with the one in $\vec{G}$,
then $c(\gamma_i)$ is the natural color of $\gamma_i$, otherwise 
$c(\gamma_i)$ is the opposite of the natural color (with understanding that the opposite
of the neutral color ``0" is ``0"). In particular, given $t\in \HB^0 \EA$ we obtain the operator
$$
D_{\gamma} \EA (t) = D_{\gamma} \EA (\underbrace{t\otimes \cdots \otimes t}_{l(\gamma)-times})=
D^{c(\gamma_{l(\gamma)})}(t)\circ \cdots \circ D^{c(\gamma_1)}(t)
$$
i.e. $t\in \HB^0 \EA$ is transformed into a composition of operators from
the space (\ref{sp-end-s}). Acting with $D_{\gamma} \EA (t)$ on $\HB^0 \EA$
has an effect of propagation of elements of $\HB^0 \EA$ through $\HT_{-w} \EA$
along the path $\gamma$.
\\
\indent
It should be noticed that the ``quantum" space $\JBBU ([e])$ also ``sees" the classical configuration $\ZE$.
The point is that the space $\HB^0 \EA$ contains distinguished set of generators
 indexed by the points $z \in \ZE$. These are obtained as follows.
Take $\delta_z$ to be the $\delta$-function on $\ZE$ supported at $z\in \ZE$,
 then its projection $\delta^0_z$ into the summand $\HB^0 \EA$ gives the aforementioned set of
generators of $\HB^0 \EA$. Thus one can say that
the ``quantum" space $\JBBU ([e])$ recovers the point-like structure of $\ZE$.
 But we obtain much more, since having
$\delta^0_z$, we can form the family of path-operators
$D_{\gamma} \EA (\delta^0_z)$ indexed by the paths of the graph $G$. This now takes a form
of ``quantum" properties of the configuration $Z_e$: a point $z \in Z_e$ represented by
$\delta^0_z$ is no longer confined to a particular location but via the action of the path-operators $D_{\gamma} \EA (\delta^0_z)$ can make its appearance elsewhere in 
$\HT_{-w} \EA$. In addition, the path-operators $D_{\gamma} \EA (\delta^0_z)$
act on $\delta^0_{z^{\prime}}$, for {\it any} $z^{\prime} \in Z_e$, and this can be regarded
as ``quantum interaction" between distinct points (particles) in $Z_e$.
\\
\indent
All these ``strange" phenomena indicate that our variety $\JUE$ of ``hidden" parameters
indeed reveals ``quantum" properties of the universal configuration $\ZDO$.
\\
\\
\indent
We turn now to a discussion of Theorem \ref{Th4}. Its essential point is to relate our
``quantum" space $\JBBU$ with the ``Higgs" space $H_{\JBBU}$ defined in (\ref{alb}) (see Theorem \ref{Th3} for some of its properties).
This relation is given by the Calabi-Yau cycle map
$CY_{\RUE}$ in (\ref{rel}). Recall that it assigns to a point $\EA \in \JBBU$ the divisor
\begin{equation}\label{CY-cycle}
\sum_{z\in Z_e} H_{z,[\AL]}
\end{equation}
where $H_{z,[\AL]}$ is a hyperplane section of $H_{\JBBU}$ intrinsically associated
to a point $(z,[\AL]) \in (Z_e,[\AL])$. In fact, the construction of $ H_{z,[\AL]}$
discussed in \S\ref{sCY} is just a geometric way to use the path-operator construction
described above. Namely, for $(z,[\AL])$ we use the path-operators
$D_{\gamma} \EA (\delta^0_z)$, where the paths $\gamma$ are the shortest (``geodesic") paths of the graph $G$
in (\ref{G}) joining the left-most (resp. the right-most) vertical edge of $G$ with the successive vertical edges
going rightwards (resp. leftwards). 
\\
\indent
Of course we could choose any other suitable collection of paths.
In fact, if we use the quantum mechanical principle of ``summing  over all histories", we
should sum over {\it all} suitable collections of paths of the graph $G(\JBBU)$.
This will lead to a ``quantum" version of the Calabi-Yau cycle map. These matters will be
considered
elsewhere (see \S6, \cite{[R1]}, for related constructions of quantum-type invariants).
\\
\\
\indent
We would like to close our discussion of string theoretic analogies by suggesting that relating the varieties $\JBBU$ and $H_{\JBBU}$ associated to
$\EE$ has a flavor of Mirror symmetry. This is based on the following observation.
\\
\indent 
  The Calabi-Yau  varieties
appearing in the construction of $CY_{\RUE}$  come in families. Namely, for a fixed configuration $Z_e$ the 
Calabi-Yau varieties $H_{z,[\AL]}\,(z\in Z_e)$ in (\ref{CY-cycle}) have moduli controlled by our
``stringy" parameter $[\AL]$ varying in the fibre $\JBBU ([e])$ of $\JBBU$ over 
$[e] \in \RUE^{\prime\prime}$.
But we have already seen that this very parameter also controls non-classical
``K\"{a}hler" structure of $Z_e$ (see the discussion following Definition \ref{ext-p}).
Thus a parameter of non-classical ``K\"{a}hler" structures on the side of
$\JBBU$ is transformed via  $CY_{\RUE}$ into a classical parameter of variation of complex structure
of hyperplane sections on the side of $H_{\JBBU}$. This can be viewed as a half of
Mirror duality between $\JBBU$ and $H_{\JBBU}$. One could speculate that the other half
of Mirror duality holds as well, i.e. the K\"{a}hler moduli of the hyperplane sections
of $H_{\JBBU}$ gets transformed into the classical, ``visible" parameter on the side
of $\JBBU$. In other words one should obtain interesting new invariants of the
configurations $Z_e$, for $[e]\in \RUE^{\prime\prime}$, from the invariants of
the families of Calabi-Yau varieties $H_{z,[\AL]}\,(z\in Z_e,\,\,[\AL] \in \JBBU ([e]))$.
Furthermore one can formulate the question of recovering the universal configuration
$\ZDO$ from the appropriate invariants of the family $H_{z,[\AL]}\,(z\in Z_e,\,\,[\AL] \in \JBBU ([e]))$. This gives an obvious analogue of the Torelli problem (see \cite{[G]}).

\section{Admissible families of clusters}
In this section we generalize our constructions to families of $0$-dimensional subschemes of $X$ or {\it clusters}\footnote{see the footnote on page
\pageref{clust}.}  
 which do not
necessarily come from vector bundles on $X$.
\\
\indent
Examining our construction of variety $\JUE$ it becomes clear that the vector bundle has been used only
to insure that every cluster of the family $p_{\UE}: \ZD \ARR \UE$ is $L$-special.
This can be formalized as follows.
\\
\indent
Fix a line bundle $\OO_X (L)$ over a smooth complex projective variety $X$ of dimension
$n\geq 2$. We assume $\OO_X (L)$ to be subject to the vanishing condition in (\ref{L-cond}).
\begin{defi}\label{clust-def}
\begin{enumerate}
\item[1)]
Let $U$ be a smooth irreducible quasi-projective variety and let
\begin{equation}\label{pU}
p_U: \ZD_U \ARR U
\end{equation}
 be a family of clusters of degree $d$ on $X$, i.e. $p_U$ is finite morphism of degree $d$
and the fibres of $p_U$ are clusters of degree $d$ on $X$. Such a family will be called a cluster
of degree $d$ on $X$ over $U$ or simply $U$-cluster of degree $d$ on $X$.
\item[2)]
A $U$-cluster of degree $d$ on $X$ is called generically smooth if the morphism 
$p_U$ in (\ref{pU}) is generically smooth, i.e. for a general $u\in U$, the fibre $Z_u$ is a set of $d$ distinct points
on $X$. Let $\stackrel{\circ}{U}$ be the locus in $U$ over which
$p_U$ is smooth. The cluster 
$$
\ZDO_U \ARR \stackrel{\circ}{U}
$$ 
is called $U$-configuration.
\item[3)]
A $U$-cluster on $X$ is called $L$-special if the fibre $Z_u$ of $p_U$ in (\ref{pU})
is $L$-special (Definition \ref{basic-def}), for every $u\in U$.
\end{enumerate}
\end{defi}
It is clear that the construction in \S2 extends to any $L$-special $U$-cluster of $X$.
\begin{pro}\label{ExtU}
Let $p_U: \ZD_U \ARR U$ be an $L$-special $U$-cluster on $X$.
There exists a torsion free sheaf $\EXTSU$ over $U$ such that for every
$u\in U$ its fibre
$\EXTSU (u) $ at $u$ is isomorphic to 
$Ext^{n-1}(\ID_{Z_u} (L),\OO_X)$, where $\ID_{Z_u}$ is the ideal sheaf of the cluster
$Z_u$.
\end{pro}
\begin{defi}
The rank of $\EXTSU$ in Proposition \ref{ExtU} is called the index of $L$-speciality
of a $U$-cluster and it will be denoted $\delta(L,U)$.
\end{defi}
As in \S2 we distinguish the Zariski open subset $U^{\prime}$ of $U$ over which
the sheaf $\EXTSU$ has constant rank $\delta(L,U)$ and define
\begin{equation}\label{jacU}
\JB_U (X;L,d) := \PP (\EXTSU \otimes \OO_{U^{\prime}})
\end{equation}
where $d$ is the degree of the $U$-cluster.
\begin{defi}
The variety $\JB_U (X;L,d)$ is called the nonabelian Jacobian of an $L$-special
$U$-cluster of degree $d$ on $X$.
\end{defi}
If no ambiguity is likely we omit some or all of the parameters $U,X,L,d$ in the above 
notation.
\\
\\
\indent
Let 
\begin{equation}\label{piU}
\pi_U: \JB_U \ARR U^{\prime}
\end{equation}
be the natural projection and let $\OO_{\JB_U} (1)$ be such that its direct image
$$
\pi_{U \ast} \OO_{\JB_U} (1) = \left(\EXTSU \otimes \OO_{U^{\prime}}\right)^{\ast}
$$
We now turn to considerations of Hodge-like structure on the sheaf
$\PIPU$. Again the only use of a vector bundle made in \S3 was the existence of 
regular extension classes (Definition \ref{ext-reg}). 
\begin{defi}\label{extU-reg}
An $L$-special $U$-cluster on $X$ is called regular if the extension group \linebreak
$Ext^{n-1}(\ID_{Z_u} (L),\OO_X)$ contains a regular extension class, for every
$u\in U$.
\end{defi}
 This condition enables us to associate a distinguished sheaf-filtration 
\begin{equation}\label{HT-filtU}
\HT_{-1} (U) \subset \ldots \subset \HT_{-k} (U)\subset \HT_{-k-1} (U) \subset\ldots
\end{equation}
of $\PIPU$. As before the filtration is torsion free and we can define the Hilbert
function $P_{L,U}$ and the weight $w(L,U)$ of the $U$-cluster (with respect to  the divisor
$L$) in exactly the same way as in 
\S3, Definition \ref{HilbUE}.
\\
\\
\indent
To have the splitting of the filtration (\ref{HT-filtU}) and hence the generalization of
 Theorems in {\bf Part II, III} of the introduction all we need is that our
$L$-special, regular $U$-cluster would be generically smooth (Definition \ref{clust-def},2)).
\begin{defi}\label{clust-adm}
An $L$-special, regular $U$-cluster on $X$ is called $L$-admissible if it is
generically smooth.
\end{defi}
\begin{cor}\label{th-cl-adm}
The Theorems \ref{Th1} -\ref{Th4} in \S0 hold for any $L$-admissible $U$-cluster on $X$.
\end{cor}

\section{Complete intersection clusters}
We will illustrate our general contructions in the case of complete intersection clusters of
sufficiently ample divisors.
\\
\indent
Let $X$ be a smooth complex projective variety of dimension $n \geq 2$ (the case $n=2$ 
has been treated in \cite{[R1]},\S5.2). We will assume that the irregularity
$q(X) = h^1 (X) =0$.
\\
\indent
Fix a very ample line bundle $\OO_X (L)$ on $X$ and consider 0-dimensional complete
intersections defined by divisors of the linear system $\mid L \mid$. Let
$Z$ be such a cluster. Then identifying $X$ with its image under the embedding defined by
$\OO_X (L)$ we obtain $Z$ by intersecting $X$ with a projective subspace of codimension
$n$ or, dually, $Z$ is defined by the global sections of $\OO_X (L)$ vanishing on $Z$, i.e.
by the subspace $H^0(\ID_{Z} (L))$ of dimension $n$ in $H^0(\OO_X (L))$.
Thus 0-dimensional complete intersections of divisors in $\mid L \mid$ are parametrized
by a Zariski open subset $U=U(L)$ of the Grassmannian  $Gr(n,H^0(\OO_X (L)))$ of 
$n$-planes in $H^0(\OO_X (L))$. If $W$ is an $n$-dimensional subspace of
$H^0(\OO_X (L))$ we denote by $[W]$ the corresponding point in $Gr(n,H^0(\OO_X (L)))$
and $Z_W$ the complete intersection subscheme of $X$ defined by $W$. In particular, it
is a well-known fact in algebraic geometry (Bertini theorem, see e.g. \cite{[H]}) that for
a general $[W] \in U$, the complete intersection $Z_W$ is smooth, i.e.
$Z_W$ is the set of $L^n$ distinct points on $X$, where $L^n$ is the self-intersection
of $L$ of degree $n$. Thus the $U$-cluster of complete intersections in $\mid L \mid$
is generically smooth (Definition \ref{clust-def}).
\\
\indent
Next we consider the question of speciality with respect to the powers of $\OO_X (L)$.
 Let $[W]\in U$ and let
$Z_W$ be the corresponding complete intersection. The ideal sheaf $\ID_{Z_W}$
of $Z_W$ admits the following Koszul resolution
\begin{equation}\label{Koz1}
\xymatrix{
  &\wedge^{2} W \otimes \OO_X (-2L)\ar[r] \ar@{{<}-}[l]& W \otimes \OO_X (-L) \ar[r]& \ID_{Z_W}\ar[r]& 0 \\
0\ar[r]& \OO_X (-nL) \ar[r]&\wedge^{n-1} W \otimes \OO_X (-(n-1)L)\ar[r]&\cdots \ar[r]&
}
\end{equation}
Let us examine the extension group
$\EXX_{Z_W} (mL) = \EXX (\ID_{Z_W} (mL), \OO_X)$, for $m>0$. By Serre duality it is
isomorphic to 
$H^1(\ID_{Z_W} (mL+K_X))^{\ast}$.
The latter group can be computed from the spectral sequence of (\ref{Koz1}) tensored with 
$\OO_X (mL+K_X)$. This yields that 
$\EXX_{Z_W} (mL)$ vanishes for all $m>n$. For $m=n$ we obtain
\begin{equation}\label{ExtnL}
\EXX_{Z_W} (nL) =H^1(\ID_{Z_W} (nL+K_X))^{\ast} = H^0(\OO_X)
\end{equation}
i.e. there exists a unique, up to a scalar multiple, nonzero extension class.
This class is of course represented by the Koszul resolution (\ref{Koz1}).
The first notrivial case from the point of view of our construction occurs for
$m=n-1$
\begin{equation}\label{Extn1L}
\EXX_{Z_W} ((n-1)L) =H^1(\ID_{Z_W} ((n-1)L+K_X))^{\ast} = H^0(\OL) /W = 
H^0(\OL) /H^0(\ID_{Z_W} (L))
\end{equation}
 Furthermore this identification implies that any section of $\OL$ whose zero-locus
is disjoint from $Z_W$ defines a regular extension class (Definition \ref{ext-reg}).
Thus we have proved the following.
\begin{pro}\label{ci1}
The $U$-cluster of complete intersections in $\mid L\mid$ is
$(n-1)L$-admissible (Definition \ref{clust-adm}) with index of speciality
$$
\delta ((n-1)L,U)=h^0 (L) - n
$$
In particular, $\delta ((n-1)L,U) \geq 2$, unless
$(X,\OL) = (\PP^n,\OO_{\PP^n} (1))$.
\end{pro}
From now on it will be assumed that the index of speciality
$\delta ((n-1)L,U) \geq 2$ and we will investigate the nature of the direct sum
decomposition in Theorem \ref{Th2}. To do this we fix $[W] \in U$ corresponding to a smooth
complete intersection $Z_W$ and a regular extension class $\AL \in \EXX_{Z_W} ((n-1)L)$.
From the identification (\ref{Extn1L}) it follows that the filtration
$\HT_{-\bullet} ([W],[\AL])$ is determined by geometry of the image of $Z_W$ under
the embedding defined by $\OL$. More precisely, the subspace
$\HT_{-k} ([W],[\AL])$ is isomorphic to the image of the obvious homomorphism
\begin{equation}\label{HT-kW}
S^k H^0(\OL) \ARR H^0(\OO_{Z_W} (kL))
\end{equation}
for every $k\geq 1$.
\\
\indent
To go further with our computations we make the following assumptions on $\OL$:
\begin{enumerate}
\item [(i)]
the line bundle $\OL$ is projectively normal, i.e. the homomorphism
\begin{equation}\label{pnor}
S^k H^0(\OL) \ARR H^0(\OO_{X} (kL))
\end{equation}
is surjective, for every $k\geq 1$;
\item [(ii)]
\begin{equation}\label{van}
H^i(\OO_{X} (kL))= 0,\,\,\,\forall k\geq 1\,\,\,and\,\,\,\forall i\geq 1
\end{equation}
\end{enumerate}
Observe that both conditions are satisfied if $L$ is a sufficiently high multiple of an arbitrary
ample divisor $D$ on $X$, i.e. $L= tD$, for an integer $t>>0$.
\\
\indent
The two assumptions above together with (\ref{HT-kW}) imply
\begin{equation}\label{PkW}
dim \HT_{-k} ([W],[\AL]) = deg Z_W - h^1(\ID_{Z_W} (kL))=L^n - h^1(\ID_{Z_W} (kL))
\end{equation}
To compute $h^1(\ID_{Z_W} (kL))$ we go back to the Koszul resolution (\ref{Koz1}) tensored
with $\OO_X (kL)$. The resulting spectral sequence together with (\ref{van}) imply that
$H^1 (\ID_{Z_W} (kL))=0, \,\forall k\geq n+1$. Thus the length of the filtration
$\HT_{-\bullet} ([W],[\AL])$ or, equivalently, the weight $w((n-1)L,U)$ of the $U$-cluster
of complete intersections in $\mid L\mid$ is at most $n+1$. Furthermore
\begin{equation}\label{PnW}
H^1 (\ID_{Z_W} (nL))= H^n (\OO_X) =H^0 (\OO_X (K_X))^{\ast}
\end{equation}
where the last equality is the Serre duality. Thus the filtration
$\HT_{-\bullet} ([W],[\AL])$ has the following form
\begin{equation}\label{HTW-filt}
\HT_{-1} ([W],[\AL]) \subset \HT_{-2} ([W],[\AL]) \subset \cdots \subset
\HT_{-n} ([W],[\AL])\subset \HT_{-n-1} ([W],[\AL]) = H^0(\OO_{Z_W})
\end{equation}
In particular, if $\AL$ is a polarizing class in $\EX_{Z_W} ((n-1)L)$, then by 
Theorem \ref{Th1} we have the following orthogonal direct sum decomposition
\begin{equation}\label{dsdW}
H^0(\OO_{Z_W}) = \bigoplus^{n}_{p=0} \HB^p ([W],[\AL])
\end{equation}
where $\HB^0 ([W],[\AL]) \cong H^0(\OL)/ W$ and $\HB^n ([W],[\AL])\cong H^0 (\OO_X (K_X))$.
Thus we proved the following.
\begin{pro}\label{ci2}
Let $X$ be a smooth complex projective variety of dimension $n\geq2$ and
$q(X) =h^1(\OO_X) =0$. 
Let
$\OL$ be a very ample line bundle on $X$ subject to (\ref{pnor}) and (\ref{van}).
Let $U=U(L)$ be the variety parametrizing 0-dimensional complete intersections
of divisors in the linear system $\mid L \mid $. Then the nonabelian Jacobian
$\JB_U (X; (n-1)L,L^n)$ is a $\PP^{\delta-1}$-bundle over $U$, where
\nolinebreak$\delta=\delta((n-1)L,U) =h^0(\OL) -n$. 
\\
\indent
Furthermore, the weight with respect to $(n-1)L$ of 
the $U$-cluster of complete intersections is \nolinebreak$w=w((n-1)L,U)\leq n+1$ and over the Zariski
open subset $\JBB_U (X; (n-1)L,L^n)$ of $\JB_U (X; (n-1)L,L^n)$ the sheaf
$\PIPU$ has the following direct sum decomposition
$$
\PIPU \otimes \OO_{\JBB_U} = \bigoplus^{n}_{p=0} \HB^p
$$
where the ranks of the first and the last summand are as follows 
$$
rk (\HB^0)= \delta((n-1)L,U) =h^0(\OL) -n,\,\,\,\,\,\,rk (\HB^n)=p_g (X)
$$
In particular, the weight $w=w((n-1)L,U)=n+1$ if and only if $p_g (X)>0$.
\end{pro}
From this result and Theorem \ref{Th3} we deduce the following.
\begin{cor}\label{ciH}
 Let $X,L$ and $U$ be as in Proposition \ref{ci2} and assume $p_g (X)>0$.
Then the nonabelian Albanese $H$ of $\JBB_U (X; (n-1)L,L^n)$ is a projective toric 
Fano variety whose 
 hyperplane sections are Calabi-Yau varieties of dimension
$n-1$.
\end{cor}
 This and the construction of the Calabi-Yau cycle map in \S5.2 imply that
behind points of a smooth complex projective variety $X$ of dimension $n\geq2$ with
irregularity $q(X) =0$ and geometric genus $p_g (X)>0$ are ``hidden" Calabi-Yau varieties of dimension
$n-1$. To ``reveal" them one needs to make a point of $X$ to be a part of
a smooth 0-dimensional complete intersection of $n$ divisors in the linear system
of sufficiently high multiple of an arbitrary ample divisor $D$.
This could be seen as a mathematical evidence of the ubiquity of 
Calabi-Yau varieties which is certainly one of the predictions of string theory.

\vspace{1cm}
\begin{flushright}
Universit\'e d'Angers\\
D\'epartement de Math\'ematiques
\\
2, boulevard Lavoisier\\
49045 ANGERS Cedex 01 \\
FRANCE\\
{\em{E-mail addres:}} reider@univ-angers.fr
\end{flushright}
  
\end{document}